\documentclass[11pt]{article}

\usepackage{amsmath,amssymb,amsfonts}
\usepackage{pstricks}
\usepackage{xr}

\numberwithin{equation}{section}
\newtheorem{thm}{Theorem}[section] 
\newtheorem{prp}[thm]{Proposition}
\newtheorem{lmm}[thm]{Lemma}  
\newtheorem{crl}[thm]{Corollary}

\def\e_ref#1{(\ref{#1})}
\def\lan{\langle}
\def\ran{\rangle}

\def\blr#1{\big\langle{#1}\big\rangle}

\def\ov#1{\overline{#1}}
\def\ti#1{\tilde{#1}}
\def\wt#1{\widetilde{#1}}

\renewcommand{\cal}{\mathcal}
\renewcommand{\frak}{\mathfrak}
\renewcommand{\Bbb}{\mathbb}

\def\A{\cal A}

\def\C{\mathbb C}
\def\cC{\cal C}
\def\E{\Bbb E}

\def\L{\Bbb L}
\def\M{\frak M}
\def\cM{\cal M}
\def\N{\cal N}
\def\O{\cal O}
\def\P{\Bbb{P}^n}
\def\bP{\Bbb{P}}
\def\Q{\Bbb Q}

\def\cZ{\mathcal Z}

\def\al{\alpha}
\def\be{\beta}
\def\io{\iota}

\def\la{\lambda}

\def\vr{\varrho}
\def\vph{\varphi}

\def\Ga{\Gamma}

\def\lra{\longrightarrow}
\def\llra{\longleftrightarrow}
\def\Lra{\Longrightarrow}

\def\eset{\emptyset}

\def\smsize#1{\begin{small}{#1}\end{small}}

\begin{document}

\title{Intersections of Tautological Classes on Blowups\linebreak
of Moduli Spaces of Genus-One Curves}

\author{Aleksey Zinger\thanks{Partially supported by an NSF Postdoctoral Fellowship}}

\date{\today}
\maketitle

\thispagestyle{empty}

\begin{abstract}
\noindent
We give two recursions for computing top intersections of tautological classes 
on blowups of moduli spaces of genus-one curves.
One of these recursions is analogous to the well-known string equation.
As shown in previous papers, these numbers are useful for computing
genus-one enumerative invariants of projective spaces and Gromov-Witten
invariants of complete intersections. 
\end{abstract}

\tableofcontents

\section{Introduction}
\label{intro_sec}

\noindent
Moduli spaces of stable curves and stable maps play a prominent role in 
algebraic geometry, symplectic topology, and string theory.
Many geometric results have been obtained by utilizing the fact that
the moduli space $\ov\M_{0,k}(\P,d)$ of degree-$d$ stable maps from
genus-zero curves with $k$ marked points into~$\P$ is a smooth unidimensional
orbi-variety of the expected dimension.
This is not the case for positive-genus moduli spaces $\ov\M_{g,k}(\P,d)$.
However, if $d\!\ge\!1$, the closure 
$$\ov\M_{1,k}^0(\P,d)\subset\ov\M_{1,k}(\P,d) $$
of the space $\M_{1,k}^0(\P,d)$ of stable maps with smooth domains  
is an irreducible orbi-variety of the expected dimension.
This component of $\ov\M_{1,k}(\P,d)$ contains all the relevant genus-one information
for the purposes of enumerative geometry and, as shown in~\cite{LZ} and~\cite{Z}, 
of the Gromov-Witten theory.\\

\noindent
For $d\!\ge\!3$, $\ov\M_{1,k}^0(\P,d)$ is singular.
A desingularization of the space~$\ov\M_{1,k}^0(\P,d)$, 
i.e.~a smooth orbi-variety $\wt\M_{1,k}^0(\P,d)$ and a~map
$$\pi\!:\wt\M_{1,k}^0(\P,d)\lra\ov\M_{1,k}^0(\P,d),$$
which is biholomorphic onto $\M_{1,k}^0(\P,d)$, is constructed in~\cite{VZ}.
Via this desingularization and the classical localization theorem of~\cite{AB},
intersections of naturally arising cohomology classes on $\ov\M_{1,k}^0(\P,d)$
can be expressed in terms of integrals of certain $\psi$-classes
on moduli spaces of genus-zero and genus-one stable curves and 
on blowups of moduli spaces of genus-one stable curves; see below for more details.
The former can be computed through two well-known recursions,
called string and dilaton equations; see Section~26.3 in~\cite{H}.
In this paper we obtain two recursions which can be used to express the latter
numbers in terms of the former ones; see Theorem~\ref{main_thm} below.
One of these recursions generalizes the genus-one string relation.\\

\noindent
If $J$ is a finite nonempty set, let $\ov\cM_{1,J}$ be the moduli space of genus-one
curves with marked points indexed by the set~$J$.
Let
$$\E\lra \ov\cM_{1,J}$$
be the Hodge line bundle of holomorphic differentials.
For each $j\!\in\!J$, we denote by
$$L_j\lra\ov\cM_{1,J}$$
the universal tangent line for the $j$th marked point and put
$$\psi_j=c_1(L_j^*) \in H^*\big(\ov\cM_{1,J};\Q\big).$$
If $(c_j)_{j\in J}$ is a tuple of integers, let
$$\blr{(c_j)_{j\in J}}_{|J|}=
\Big\lan\prod_{j\in J}\psi_j^{c_j},\ov\cM_{1,J}\Big\ran.$$\\

\noindent
Let $I$ and $J$ be two finite sets, not both empty.
The inductive procedure of Subsection~\ref{g1desing-curve1bl_subs} in~\cite{VZ}, 
which is reviewed in Subsection~\ref{blowup_subs} below, constructs a blowup
$$\pi\!: \wt\cM_{1,(I,J)}\lra \ov\cM_{1,I\sqcup J}$$
of $\ov\cM_{1,I\sqcup J}$ along natural subvarieties and their proper transforms.
In addition, it describes $|I|\!+\!1$ line bundles 
$$\ti\E,\, \ti{L}_i\lra \wt\cM_{1,(I,J)}, \qquad i\!\in\!I,$$
and $|I|$ nowhere vanishing sections 
$$\ti{s}_i\in\Ga\big(\wt\cM_{1,(I,J)};\ti{L}_i^*\!\otimes\!\ti\E^*\big),
\qquad i\!\in\!I.$$
These line bundles are obtained by twisting $\E$ and~$L_i$.
Since the sections $\ti{s}_i$ do not vanish,
all $|I|\!+\!1$ line bundles $\ti{L}_i$ and $\ti\E^*$ are  explicitly isomorphic.
They will be denoted by
$$\L\lra \wt\cM_{1,(I,J)}$$ 
and called the universal tangent bundle.
Let
$$\ti\psi=c_1(\L^*)\in H^2\big(\wt\cM_{1,(I,J)};\Q\big)$$
be the corresponding ``$\psi$-class" on $\wt\cM_{1,(I,J)}$.
If $(\ti{c},(c_j)_{j\in J})$ is a tuple of integers, we put
\begin{equation}\label{psinums_e2}
\blr{\ti{c};(c_j)_{j\in J}}_{(|I|,|J|)}=
\Big\lan \ti\psi^{\ti{c}}
\cdot\prod_{j\in J}\pi^*\psi_j^{c_j},\wt\cM_{1,(I,J)}\Big\ran.
\end{equation}
If $\sum_{j\in J}\!c_j\!\neq\!|J|$, $\ti{c}\!<\!0$, or $c_j\!<\!0$ for some $j\!\in\!J$,
we define this number to be zero.\\

\begin{thm}
\label{main_thm} 
Suppose  $I$ and $J$ are finite sets, such that $|I|\!+\!|J|\!\ge\!2$, and 
$(\ti{c},(c_j)_{j\in J})$ is a tuple of integers.
If $c_{j^*}\!=\!0$ for some $j^*\!\in\!J$,
$$\blr{\ti{c};(c_j)_{j\in J}}_{(|I|,|J|)}= 
|I|\blr{\ti{c}\!-\!1;(c_j)_{j\in J-\{j^*\}}}_{(|I|,|J|-1)}
+\!\sum_{j\in J-\{j^*\}}\!\!\!\!\!\!
      \blr{\ti{c};c_j\!-\!1,(c_{j'})_{j'\in J-\{j^*,j\}}}_{(|I|,|J|-1)}.$$
If $I\!\neq\!\eset$ and $c_j\!>\!0$ for all $j\!\in\!J$,
$$\blr{\ti{c};(c_j)_{j\in J}}_{(|I|,|J|)}= 
\blr{\ti{c};(c_j)_{j\in J}}_{(|I|-1,|J|+1)}.$$
\end{thm}

\begin{crl}
\label{psiclass_crl}
If $I$ and $J$ are finite sets and $I\!\neq\!\eset$, then
$$\blr{\ti\psi^{|I|+|J|},\wt\cM_{I,J}}
=\frac{1}{24}\cdot |I|^{|J|} \cdot (|I|\!-\!1)!$$\\
\end{crl}

\noindent
We recall that 
$$\blr{\psi,\ov\cM_{1,1}}=\frac{1}{24}.$$
Thus, Corollary~\ref{psiclass_crl} follows from Theorem~\ref{main_thm} by applying 
the first recursion $|J|$ times and then 
the second recursion followed by the first $|I|\!-\!1$ times.\\

\noindent
It is immediate from the construction of Subsection~\ref{blowup_subs} below that
$$I=\eset  \qquad\Lra\qquad \wt\cM_{1,(I,J)}=\ov\cM_{1,I\sqcup J}
\quad\hbox{and}\quad \ti\psi= \la\!\equiv\!c_1(\E).$$
Thus, the two recursions of Theorem~\ref{main_thm}, along with 
the string and dilaton equations, provide a straightforward algorithm 
for computing all numbers~\e_ref{psinums_e2}.\\

\noindent
Note that if we take $\ti{c}\!=\!0$ in the first equation of Theorem~\ref{main_thm},
we recover the genus-one string recursion for $\psi$-classes.
This equation is proved in Subsection~\ref{pfoutline_subs}
by an argument similar to the standard proof of the string recursion.
In particular, we consider the forgetful morphism
$$f\!:\ov\cM_{1,I\sqcup J}\lra  \ov\cM_{1,I\sqcup(J-\{j^*\})}.$$
We show in Subsection~\ref{lift_subs} that it lifts to a morphism on the blowups,
$$\ti{f}\!:\wt\cM_{1,(I,J)}\lra  \wt\cM_{1,(I,J-\{j^*\})}.$$
The first recursion of Theorem~\ref{main_thm} is obtained by comparing
$\ti\psi$ with $\ti{f}^*\ti\psi$.
On the other hand, the second equation of Theorem~\ref{main_thm} follows easily 
from the relevant definitions, which are reviewed in Subsection~\ref{blowup_subs}.
The reason is that the blowups of $\ov\cM_{1,I\sqcup J}$ corresponding to 
the two sides of this equation differ by blowups along loci on which
$\prod_{j\in J}\psi_j$ vanishes; see the end of Subsection~\ref{blowup_subs}.\\

\noindent
In~\cite{VZ},\ it is observed that if $n\!<\!m$, the natural embedding
$$\io\!:\ov\M_{1,k}(\P,d) \lra\ov\M_{1,k}(\bP^m,d)$$
induced by the inclusion $\P\!\lra\!\bP^m$ lifts to an embedding on
the desingularizations:
$$\ti\io\!:\wt\M_{1,k}^0(\P,d)\lra \wt\M_{1,k}^0(\bP^m,d).$$
Proceeding analogously to Subsection~\ref{lift_subs}, one can show that if $k\!>\!0$
the forgetful morphism
$$f\!:\ov\M_{1,k}(\P,d) \lra \ov\M_{1,k-1}(\P,d)$$
lifts to a morphism
$$\ti{f}\!:\wt\M_{1,k}^0(\P,d) \lra \wt\M_{1,k-1}^0(\P,d).$$
Thus, the desingularization $\wt\M_{1,k}^0(\P,d)$ of $\ov\M_{1,k}^0(\P,d)$
constructed in~\cite{VZ} respects at least two properties
that play a central role in the Gromov-Witten theory; see Figure~\ref{lift_fig0}.

\begin{figure}
\begin{pspicture}(-3,-1.8)(10,1.3)
\psset{unit=.4cm}
\rput(0,0){\smsize{$\wt\M_{1,k}^0(\P,d)$}}
\rput(9,0){\smsize{$\wt\M_{1,k}^0(\bP^m,d)$}}
\rput(0,-4){\smsize{$\ov\M_{1,k}(\P,d)$}}
\rput(9,-4){\smsize{$\ov\M_{1,k}(\bP^m,d)$}}
\psline[linestyle=dashed]{->}(2.5,0)(6.5,0)
\psline{->}(2.5,-4)(6.5,-4)\psline{->}(0,-.8)(0,-3.2)\psline{->}(9,-.8)(9,-3.2)
\rput(4.5,.5){\smsize{$\ti\io$}}\rput(4.5,-3.5){\smsize{$\io$}}
\rput(-.5,-2){\smsize{$\pi$}}\rput(8.5,-2){\smsize{$\pi$}}
\rput(18,0){\smsize{$\wt\M_{1,k}^0(\P,d)$}}
\rput(27,0){\smsize{$\wt\M_{1,k-1}^0(\P,d)$}}
\rput(18,-4){\smsize{$\ov\M_{1,k}(\P,d)$}}
\rput(27,-4){\smsize{$\ov\M_{1,k-1}(\P,d)$}}
\psline[linestyle=dashed]{->}(20.5,0)(24,0)
\psline{->}(20.5,-4)(24,-4)\psline{->}(18,-.8)(18,-3.2)\psline{->}(27,-.8)(27,-3.2)
\rput(22.5,.6){\smsize{$\ti{f}$}}\rput(22.5,-3.4){\smsize{$f$}}
\rput(17.5,-2){\smsize{$\pi$}}\rput(26.5,-2){\smsize{$\pi$}}
\end{pspicture}
\caption{Some Natural Properties of $\wt\M_{1,k}^0(\P,d)$}
\label{lift_fig0}
\end{figure}

\section{Preliminaries}
\label{prelim_subs}

\subsection{Blowup Construction}
\label{blowup_subs}

\noindent
If $I$ is a finite set, let
\begin{equation}\label{g0and1curvsubv_e}
\A_1(I) =\big\{\big(I_P,\{I_k\!:k\!\in\!K\}\big)\!: 
K\!\neq\!\eset;~I\!=\!\bigsqcup_{k\in\{P\}\sqcup K}\!\!\!\!\!\!\!I_k;~
|I_k|\!\ge\!2 ~\forall\, k\!\in\!K\big\}.
\end{equation}
Here $P$ stands for ``principal" (component).
If $\rho\!=\!(I_P,\{I_k\!:k\!\in\!K\})$ is an element of $\A_1(I)$, 
we denote by $\cM_{1,\rho}$ the subset of $\ov\cM_{1,I}$
consisting of the stable curves~$\cC$ such~that\\
${}\quad$ (i) $\cC$ is a union of a smooth torus and $|K|$ projective lines,
indexed by~$K$;\\
${}\quad$ (ii) each line is attached directly to the torus;\\
${}\quad$ (iii) for each $k\!\in\!K$,
the marked points on the line corresponding to $k$ are indexed by $I_k$.\\
Let $\ov\cM_{1,\rho}$ be the closure of $\cM_{1,\rho}$ in~$\ov\cM_{1,I}$.
Figure~\ref{g1curv_fig} illustrates this definition,
from the points of view of symplectic topology and of algebraic geometry.
In the first diagram, each circle represents a sphere, or~$\bP^1$.
In the second diagram, the irreducible components of $\cC$ are represented by curves,
and the integer next to each component shows its genus. 
It is well-known that each space $\ov\cM_{1,\rho}$ is 
a smooth subvariety of~$\ov\cM_{1,I}$.\\

\begin{figure}
\begin{pspicture}(-1.1,-1.8)(10,1.3)
\psset{unit=.4cm}
\rput{45}(0,-4){\psellipse(5,-1.5)(2.5,1.5)\pscircle*(7.5,-1.5){.2}\pscircle*(2.5,-1.5){.2}
\psarc[linewidth=.05](5,-3.3){2}{60}{120}\psarc[linewidth=.05](5,0.3){2}{240}{300}
\pscircle(5,-4){1}\pscircle*(5,-3){.15}
\pscircle*(5,-5){.2}\pscircle*(4,-4){.2}\pscircle*(6,-4){.2}
\pscircle(6.83,.65){1}\pscircle*(6.44,-.28){.15}
\pscircle*(5.9,1.04){.2}\pscircle*(7.76,.26){.2}
\pscircle(3.17,.65){1}\pscircle*(3.56,-.28){.15}
\pscircle*(4.1,1.04){.2}\pscircle*(2.24,.26){.2}}
\rput(7,.4){\smsize{$i_1$}}\rput(2.4,-3.7){\smsize{$i_2$}}
\rput(1.1,-2.7){\smsize{$i_3$}}\rput(2,.3){\smsize{$i_4$}}
\rput(3.1,.5){\smsize{$i_5$}}\rput(6,1.8){\smsize{$i_6$}}
\rput(7.7,-2.5){\smsize{$i_7$}}\rput(7.7,-4.2){\smsize{$i_8$}}
\rput(5.3,-4.5){\smsize{$i_9$}}
\psarc(15,-1){3}{-60}{60}\pscircle*(16.93,1.3){.2}\pscircle*(16.93,-3.3){.2}
\rput(17.6,1.4){\smsize{$i_1$}}
\rput(17.6,-3.4){\smsize{$i_2$}}
\psline(16.8,0)(22.05,1.25)\pscircle*(18.9,.5){.2}\pscircle*(21,1){.2}
\rput(18.9,1.2){\smsize{$i_3$}}\rput(21,1.7){\smsize{$i_4$}}
\psline(17,-1)(22,-1)\pscircle*(19.5,-1){.2}\pscircle*(21,-1){.2}
\rput(19.5,-.3){\smsize{$i_5$}}\rput(21,-.3){\smsize{$i_6$}}
\psline(16.8,-2)(22.05,-3.25)\pscircle*(18.9,-2.5){.2}
\pscircle*(19.95,-2.75){.2}\pscircle*(21,-3){.2}
\rput(18.9,-1.8){\smsize{$i_7$}}\rput(19.95,-2.05){\smsize{$i_8$}}
\rput(21,-2.3){\smsize{$i_9$}}
\rput(16.1,-3.5){$1$}\rput(22.5,1.2){$0$}\rput(22.4,-1){$0$}\rput(22.5,-3.2){$0$}
\rput(32,-1){\begin{small}
\begin{tabular}{l}
$I_P\!=\!\{i_1,i_2\}$\\
$K\!=\!\{1,2,3\}$\\
$I_1\!=\!\{i_3,i_4\}$\\
$I_2\!=\!\{i_5,i_6\}$\\
$I_3\!=\!\{i_7,i_8,i_9\}$\\
\end{tabular}\end{small}}
\end{pspicture}
\caption{A Typical Element of $\ov\cM_{1,\rho}$}
\label{g1curv_fig}
\end{figure}

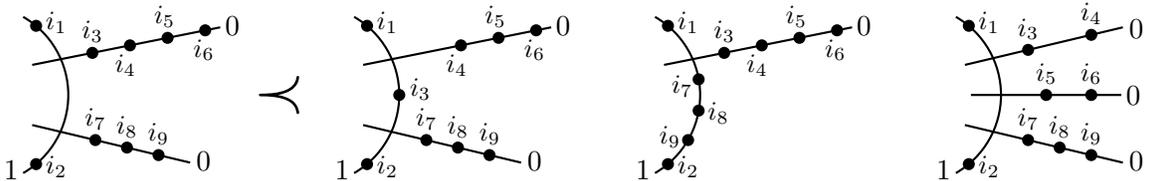
\begin{figure}
\begin{pspicture}(-1.1,-1.8)(10,1.3)
\psset{unit=.4cm}
\psarc(-2,-1){3}{-60}{60}\pscircle*(-.07,1.3){.2}\pscircle*(-.07,-3.3){.2}
\rput(.6,1.4){\begin{small}$i_1$\end{small}}
\rput(.6,-3.4){\begin{small}$i_2$\end{small}}
\psline(-.2,0)(6.05,1.25)\pscircle*(1.8,.4){.2}\pscircle*(3.05,.65){.2}
\pscircle*(4.3,.9){.2}\pscircle*(5.55,1.15){.2}
\rput(1.8,1.1){\begin{small}$i_3$\end{small}}
\rput(2.9,0){\begin{small}$i_4$\end{small}}
\rput(4.2,1.6){\begin{small}$i_5$\end{small}}
\rput(5.5,.5){\begin{small}$i_6$\end{small}}
\psline(-.2,-2)(5.05,-3.25)\pscircle*(1.9,-2.5){.2}
\pscircle*(2.95,-2.75){.2}\pscircle*(4,-3){.2}
\rput(1.9,-1.8){\begin{small}$i_7$\end{small}}
\rput(2.95,-2.05){\begin{small}$i_8$\end{small}}
\rput(4,-2.3){\begin{small}$i_9$\end{small}}
\rput(-.9,-3.5){$1$}\rput(6.45,1.3){$0$}\rput(5.5,-3.2){$0$}
\rput(8,-1){\begin{Huge}$\prec$\end{Huge}}
\psarc(9,-1){3}{-60}{60}\pscircle*(10.93,1.3){.2}\pscircle*(10.93,-3.3){.2}
\rput(11.6,1.4){\begin{small}$i_1$\end{small}}
\rput(11.6,-3.4){\begin{small}$i_2$\end{small}}
\psline(10.8,0)(17.05,1.25)\pscircle*(14.05,.65){.2}
\pscircle*(15.3,.9){.2}\pscircle*(16.55,1.15){.2}
\pscircle*(12,-1){.2}\rput(12.7,-.8){\begin{small}$i_3$\end{small}}
\rput(13.9,0){\begin{small}$i_4$\end{small}}
\rput(15.2,1.6){\begin{small}$i_5$\end{small}}
\rput(16.5,.5){\begin{small}$i_6$\end{small}}
\psline(10.8,-2)(16.05,-3.25)\pscircle*(12.9,-2.5){.2}
\pscircle*(13.95,-2.75){.2}\pscircle*(15,-3){.2}
\rput(12.9,-1.8){\begin{small}$i_7$\end{small}}
\rput(13.95,-2.05){\begin{small}$i_8$\end{small}}
\rput(15,-2.3){\begin{small}$i_9$\end{small}}
\rput(10.1,-3.5){$1$}\rput(17.45,1.3){$0$}\rput(16.5,-3.2){$0$}
\psarc(19,-1){3}{-60}{60}\pscircle*(20.93,1.3){.2}\pscircle*(20.93,-3.3){.2}
\rput(21.6,1.4){\begin{small}$i_1$\end{small}}
\rput(21.6,-3.4){\begin{small}$i_2$\end{small}}
\psline(20.8,0)(27.05,1.25)\pscircle*(22.8,.4){.2}\pscircle*(24.05,.65){.2}
\pscircle*(25.3,.9){.2}\pscircle*(26.55,1.15){.2}
\rput(22.8,1.1){\begin{small}$i_3$\end{small}}
\rput(23.9,0){\begin{small}$i_4$\end{small}}
\rput(25.2,1.6){\begin{small}$i_5$\end{small}}
\rput(26.5,.5){\begin{small}$i_6$\end{small}}
\pscircle*(21.95,-.48){.2}\rput(21.4,-.7){\begin{small}$i_7$\end{small}}
\pscircle*(21.95,-1.52){.2}\rput(22.6,-1.6){\begin{small}$i_8$\end{small}}
\pscircle*(21.6,-2.5){.2}\rput(21,-2.4){\begin{small}$i_9$\end{small}}
\rput(20.1,-3.5){$1$}\rput(27.45,1.3){$0$}
\psarc(29,-1){3}{-60}{60}\pscircle*(30.93,1.3){.2}\pscircle*(30.93,-3.3){.2}
\rput(31.6,1.4){\begin{small}$i_1$\end{small}}
\rput(31.6,-3.4){\begin{small}$i_2$\end{small}}
\psline(30.8,0)(36.05,1.25)\pscircle*(32.9,.5){.2}\pscircle*(35,1){.2}
\rput(32.9,1.2){\begin{small}$i_3$\end{small}}
\rput(35,1.7){\begin{small}$i_4$\end{small}}
\psline(31,-1)(36,-1)\pscircle*(33.5,-1){.2}\pscircle*(35,-1){.2}
\rput(33.5,-.3){\begin{small}$i_5$\end{small}}
\rput(35,-.3){\begin{small}$i_6$\end{small}}
\psline(30.8,-2)(36.05,-3.25)\pscircle*(32.9,-2.5){.2}
\pscircle*(33.95,-2.75){.2}\pscircle*(35,-3){.2}
\rput(32.9,-1.8){\begin{small}$i_7$\end{small}}
\rput(33.95,-2.05){\begin{small}$i_8$\end{small}}
\rput(35,-2.3){\begin{small}$i_9$\end{small}}
\rput(30.1,-3.5){$1$}\rput(36.5,1.2){$0$}\rput(36.4,-1){$0$}\rput(36.5,-3.2){$0$}
\end{pspicture}
\caption{Examples of Partial Ordering~\e_ref{partorder_e}}
\label{partorder_fig}
\end{figure}

\noindent
We define a partial ordering on the set $\A_1(I)\!\sqcup\!\{(I,\eset)\}$  by setting
\begin{equation}\label{partorder_e}
\rho'\!\equiv\!\big(I_P',\{I_k'\!: k\!\in\!K'\}\big)
\prec \rho\!\equiv\!\big(I_P,\{I_k\!: k\!\in\!K\}\big)
\end{equation}
if $\rho'\!\neq\!\rho$ and there exists a map $\vph\!:K\!\lra\!K'$
such that $I_k\!\subset\!I_{\vph(k)}'$ for all $k\!\in\!K$.
This condition means that the elements of $\cM_{1,\rho'}$ can be obtained
from the elements of $\cM_{1,\rho}$ by moving more points onto the bubble components
or combining the bubble components; see Figure~\ref{partorder_fig}.\\

\noindent
Let $I$ and $J$ be finite sets such that $I$ is not empty and $|I|\!+\!|J|\!\ge\!2$. 
We put
\begin{equation*}
\A_1(I,J) =\big\{\big((I_P\!\sqcup\!J_P),\{I_k\!\sqcup\!J_k\!: k\!\in\!K\}\big)\!\in\!\A_1(I\!\sqcup\!J)\!:
~I_k\!\neq\!\eset~ \forall\, k\!\in\!K \big\}.
\end{equation*}
We note that if $\vr\!\in\!\A_1(I\!\sqcup\!J)$, then $\vr\!\in\!\A_1(I,J)$ 
if and only if every bubble component of an element of $\cM_{1,\vr}$ 
carries at least one element of~$I$.
The partially ordered set $(\A_1(I,J),\prec)$ has a unique minimal element
$$\vr_{\min}\equiv\big(\eset,\{I\!\sqcup\!J\}\big).$$
Let $<$ be an ordering on $\A_1(I,J)$ extending the partial ordering $\prec$.
We denote the corresponding maximal element by~$\vr_{\max}$.
If $\vr\!\in\!\A_1(I,J)$, we~put
\begin{equation}\label{minusdfn_e}
\vr\!-\!1=
\begin{cases}
\max\{\vr'\!\in\!\A_1(I,J)\!: \vr'\!<\!\vr\},&
\hbox{if}~ \vr\!\neq\!\vr_{\min};\\
0,& \hbox{if}~\vr\!=\!\vr_{\min},
\end{cases}
\end{equation}
where the maximum is taken with respect to the ordering $<$.\\

\noindent
The starting data for the blowup construction of 
Subsection~\ref{g1desing-curve1bl_subs} in~\cite{VZ} is given~by
\begin{gather*}
\ov\cM_{1,(I,J)}^0=\ov\cM_{1,I\sqcup J},\qquad
\ov\cM_{1,\vr}^0=\ov\cM_{1,\vr} \quad \forall\,\vr\!\in\!\A_1(I,J),\\
\E_0=\E\lra \ov\cM_{1,(I,J)}^0, \qquad\hbox{and}\qquad
L_{0,i}\!=\!L_i\lra\ov\cM_{1,(I,J)}^0 \quad\forall\,i\!\in\!I.
\end{gather*}
Suppose $\vr\!\in\!\A_1(I,J)$ and we have constructed\\
${}\quad$ ($I1$) a blowup $\pi_{\vr-1}\!:\ov\cM_{1,(I,J)}^{\vr-1}\!\lra\!\ov\cM_{1,(I,J)}^0$
of $\ov\cM_{1,(I,J)}^0$ such that $\pi_{\vr-1}$ is one-to-one outside of the preimages
of the spaces $\ov\cM_{1,\vr'}^0$ with $\vr'\!\le\!\vr-1$;\\
${}\quad$ ($I2$) line bundles $L_{\vr-1,i}\!\lra\!\ov\cM_{1,(I,J)}^{\vr-1}$ 
for $i\!\in\!I$ and $\E_{\vr-1}\!\lra\!\ov\cM_{1,(I,J)}^{\vr-1}$.\\
For each $\vr^*\!>\!\vr\!-\!1$, let $\ov\cM_{1,\vr^*}^{\vr-1}$ 
be the proper transform of~$\ov\cM_{1,\vr^*}^0$ in~$\ov\cM_{1,(I,J)}^{\vr-1}$.\\

\noindent
If $\vr\!\in\!\A_1(I,J)$ is as above, let 
$$\ti\pi_{\vr}\!:\ov\cM_{1,(I,J)}^{\vr}\lra\ov\cM_{1,(I,J)}^{\vr-1}$$ 
be the blowup of $\ov\cM_{1,(I,J)}^{\vr-1}$ along $\ov\cM_{1,\vr}^{\vr-1}$.
We denote by $\ov\cM^{\vr}_{1,\vr}$ the corresponding exceptional divisor.
If $\vr^*\!>\!\vr$, let $\ov\cM_{1,\vr^*}^{\vr}\!\subset\!\ov\cM_{1,(I,J)}^{\vr}$ be 
the proper transform of~$\ov\cM_{1,\vr^*}^{\vr-1}$.
If 
$$\vr=\big((I_P\!\sqcup\!J_P),\{I_k\!\sqcup\!J_k\!: k\!\in\!K\}\big)\!\in\!\A_1(I\!\sqcup\!J)
\qquad\hbox{and}\qquad i\!\in\!I,$$ 
we~put
\begin{equation}\label{bundtwist_e}
L_{\vr,i}=\begin{cases}
\ti\pi_{\vr}^*L_{\vr-1,i},& \hbox{if}~ i\!\not\in\!I_P;\\
\ti\pi_{\vr}^*L_{\vr-1,i}\otimes\O(-\ov\cM^{\vr}_{1,\vr}),& \hbox{if}~ i\!\in\!I_P;
\end{cases}  \qquad
\E_{\vr}=\ti\pi_{\vr}^*\,\E_{\vr-1}\otimes\O(\ov\cM^{\vr}_{1,\vr}).
\end{equation}
It is immediate that the requirements ($I1$) and ($I2$),
with $\vr\!-\!1$ replaced by~$\vr$, are satisfied.\\

\noindent
We conclude the blowup construction after $|\vr_{\max}|$ steps.
Let
$$\wt\cM_{1,(I,J)}=\ov\cM_{1,(I,J)}^{\vr_{\max}}, \qquad
\ti{L}_i=L_{\vr_{\max},i} ~~\forall\, i\!\in\!I, \qquad
\ti{\E}=\E_{\vr_{\max}}.$$
By Lemma~\ref{g1desing-curve1bl_lmm} in~\cite{VZ}, 
the end result of this blowup construction 
is well-defined, i.e.~independent of the choice of an ordering $<$ extending 
the partial ordering~$\prec$.
The reason is that different extensions of the partial order~$\prec$ correspond
to different orders of blowups along disjoint subvarieties.
By the inductive assumption~($I4$) in Subsection~\ref{g1desing-curve1bl_subs} of~\cite{VZ},
there is a natural isomorphism between the line bundles $\ti{L}_i$ and~$\ti\E^*$.
Thus, these line bundles are the same. We denote them by~$\L$.\\

\noindent
We are now ready to verify the second equation in Theorem~\ref{main_thm}.
If $i^*\!\in\!I$,
\begin{gather*}
\A_1\big(I\!-\!\{i^*\},J\!\sqcup\!\{i^*\}\big)  \subset \A_1(I,J)
\qquad\hbox{and}\\
\A_1(I,J)\!-\!\A_1\big(I\!-\!\{i^*\},J\!\sqcup\!\{i^*\}\big)
=\big\{\!\vr\!=\!\big(I_P\!\sqcup\!J_P,
\big\{\{i^*\}\!\sqcup\!J_1\big\}\!\sqcup\!\{I_k\!\sqcup\!J_k\!:k\!\in\!K'\}\big)
\!\in\!\A_1(I\!\sqcup\!J)\big\}.
\end{gather*}
With $\vr$ as above, we have a natural isomorphism
$$\ov\cM_{1,\vr}\approx \ov\cM_{1,\bar\vr}\times\ov\cM_{0,\{q,i^*\} \sqcup J_1},
\qquad\hbox{where} \quad
\bar\vr= \big(I_P\!\sqcup\!J_P\!\sqcup\!\{p\},\{I_k\!\sqcup\!J_k\!:k\!\in\!K'\}\big).$$
Let 
$$\pi_2\!: \ov\cM_{1,\vr}\lra\ov\cM_{0,\{q,i^*\}\sqcup J_1}$$
be the projection map.
By definition, 
$$\psi_j\big|_{\ov\cM_{1,\vr}}=\pi_2^*\psi_j \qquad \forall\, j\!\in\!J_1
\qquad\Lra\qquad
\prod_{j\in J_1}\!\!\psi_j\,\big|_{\ov\cM_{1,\vr}}=
\pi_2^*\prod_{j\in J_1}\!\!\psi_j=\pi_2^*0=0,$$
since the dimension of $\ov\cM_{0,\{q,i^*\}\sqcup J_1}$ is $|J_1|\!-\!1$.
It follows that
$$\prod_{j\in J}\!\psi_j \,\big|_{\ov\cM_{1,\vr}}=0 \qquad
\forall\,\vr\!\in\!\A_1(I,J)\!-\!\A_1\big(I\!-\!\{i^*\},J\!\sqcup\!\{i^*\}\big).$$
Thus, the constructions of $\ti\psi\!\equiv\!c_1(\ti\E)$ from $\la\!\equiv\!c_1(\E_0)$
for $\wt\cM_{1,(I-\{i^*\},J\sqcup\{i^*\})}$ and $\wt\cM_{1,(I,J)}$
differ by varieties along which $\prod_{j\in J}\!\psi_j^{c_j}$ vanishes, 
as long as $c_j\!>\!0$ for all $j\!\in\!J$.
We conclude that
$$\Big\lan \ti\psi^{\ti{c}}
\cdot\prod_{j\in J}\pi^*\psi_j^{c_j},\wt\cM_{1,(I,J)}\Big\ran
=\Big\lan \ti\psi^{\ti{c}}
\cdot\prod_{j\in J}\pi^*\psi_j^{c_j},\wt\cM_{1,(I-\{i^*\},J\sqcup\{i^*\})}\Big\ran$$
whenever $c_j\!>\!0$ for all $j\!\in\!J$, as needed.

\subsection{Outline of Proof of First Recursion in Theorem~\ref{main_thm}}
\label{pfoutline_subs}

\noindent
In this subsection we state three results, Proposition~\ref{lift_prp}
and Lemmas~\ref{psirestr_lmm1} and~\ref{psirestr_lmm2}, that lead 
in a straightforward way to the first recursion of Theorem~\ref{main_thm}.
They are proved in the next section.\\

\noindent
If $I$ is a finite set and $i,j$ are distinct elements of $I$, let
$$\rho_{ij}=\big(I\!-\!\{i,j\},\{\{i,j\}\}\big) \in \A_1(I).$$
There is a natural decomposition
\begin{equation}\label{decomp_e}
\ov\cM_{1,\vr_{ij}} = \ov\cM_{1,(I-\{i,j\})\sqcup\{p\}}\times\ov\cM_{0,\{q,i,j\}}.
\end{equation}
The second component is a one-point space.
Let
\begin{equation}\label{projmap_e}
\pi_P,\pi_B\!: \ov\cM_{1,\vr_{ij}} \lra 
\ov\cM_{1,(I-\{i,j\})\sqcup\{p\}},  \ov\cM_{0,\{q,i,j\}}
\end{equation}
be the two projection maps.
Here $P$ and $B$ stand for ``principal" and ``bubble" (components).
It is immediate that
\begin{gather}
\label{larestr_e1}
\la|_{\ov\cM_{1,\vr_{ij}}}=\pi_P^*\la \qquad\hbox{and}\\
\label{psirestr_e1}
\psi_{j'}\big|_{\ov\cM_{1,\vr_{ij}}}=\begin{cases}
\pi_P^*\psi_{j'},&\hbox{if}~j'\!\neq\!i,j;\\
\pi_B^*\psi_{j'}\!=\!0,&\hbox{if}~j'\!=\!i,j;
\end{cases} \qquad\forall\,j'\!\in\!I.
\end{gather}
In the $j'\!=\!i,j$ case the restriction of $\psi_{j'}$ 
vanishes because the second component is zero-dimensional.\\

\noindent
If $I$ is a finite set, $|I|\!\ge\!2$, and $j^*\!\in\!I$, there is a natural forgetful
morphism
$$f\!: \ov\cM_{1,I}\lra \ov\cM_{1,I-\{j^*\}}.$$
It is obtained by dropping the marked point $j^*$ from every element of $\ov\cM_{1,I}$
and contracting the unstable components of the resulting curve.
It is straightforward to check that
\begin{gather}\label{lapullback_e1}
\la=f^*\la \qquad\hbox{and}\\
\label{psipullback_e1}
\psi_j=f^*\psi_j+\ov\cM_{1,\vr_{jj^*}} \quad\Lra\quad
f^*\psi_j|_{\ov\cM_{1,\vr_{jj^*}}}=\pi_P^*\psi_p
\qquad\forall\,j\!\in\!I\!-\!\{j^*\};
\end{gather}
see Chapter~25 in \cite{H}, for example.
From \e_ref{psirestr_e1} and \e_ref{psipullback_e1}, we find that
\begin{equation}\label{psipullback_e1b}
\psi_j^{c_j}=\psi_j^{c_j-1}\big(f^*\psi_j+\ov\cM_{1,\vr_{jj^*}}\big)
=f^*\psi_j^{c_j}+\big(\pi_P^*\psi_p^{c_j-1}\big)\!\cap\!\ov\cM_{1,\vr_{jj^*}}
\quad\forall\,j\!\in\!I\!-\!\{j^*\},\,c_j>0.
\end{equation}\\

\noindent
If $I$ and $J$ are finite sets, $i\!\in\!I$, and $j\!\in\!J$,
then $\ov\cM_{1,\vr_{ij}}$ is a divisor in $\ov\cM_{1,I\sqcup J}$.
Thus, in the notation of the previous subsection,
$$\ov\cM_{1,\vr_{ij}}^{\vr_{ij}} = \ov\cM_{1,\vr_{ij}}^{\vr_{ij}-1}.$$
Since $\vr_{ij}$ is a maximal element of $(A_1(I,J),\prec)$,
the blowup loci at the stages of the construction described in 
Subsection~\ref{blowup_subs} that follow the blowup along 
$\ov\cM_{1,\vr_{ij}}^{\vr_{ij}-1}$ are disjoint from $\ov\cM_{1,\vr_{ij}}^{\vr_{ij}}$.
Thus, we can view $\ov\cM_{1,\vr_{ij}}^{\vr_{ij}}$ as a divisor in
$\wt\cM_{1,(I,J)}$. We denote it by $\wt\cM_{1,\vr_{ij}}$.
If $i,j\!\in\!J$, $\ov\cM_{1,\vr_{ij}}$ is also a divisor in $\ov\cM_{1,I\sqcup J}$.
Thus, its proper transform $\ov\cM_{1,\vr_{ij}}^{\vr}$ in $\ov\cM_{1,(I,J)}^{\vr}$
is a divisor for every $\vr\!\in\!\A_1(I,J)$. 
Let
$$\wt\cM_{1,\vr_{ij}}=\ov\cM_{1,\vr_{ij}}^{\vr_{\max}}\subset\wt\cM_{1,(I,J)}.$$

\begin{prp}
\label{lift_prp}
Suppose $I$ and $J$ are finite sets such that $|I|\!+\!|J|\!\ge\!2$ and $j^*\!\in\!J$.
If  $$\pi\!: \wt\cM_{1,(I,J)}\lra\ov\cM_{1,I\sqcup J} \qquad\hbox{and}\qquad
\pi\!: \wt\cM_{1,(I,J-\{j^*\})}\lra\ov\cM_{1,I\sqcup (J-\{j^*\})}$$
are blowups as in Subsection~\ref{blowup_subs}, the forgetful map
$$f\!: \ov\cM_{1,I\sqcup J}\lra \ov\cM_{1,I\sqcup (J-\{j^*\})}$$
lifts to a morphism 
$$\ti{f}\!: \wt\cM_{1,(I,J)}\lra \wt\cM_{1,(I,J-\{j^*\})};$$
see Figure~\ref{lift_fig1}.
Furthermore,
$$\ti\psi=\ti{f}^*\ti\psi+\sum_{i\in I}\wt\cM_{1,\vr_{ij^*}}.$$
\end{prp}

\begin{figure}
\begin{pspicture}(-6,-1.8)(10,1.3)
\psset{unit=.4cm}
\rput(0,0){\smsize{$\wt\cM_{1,(I,J)}$}}
\rput(9,0){\smsize{$\wt\cM_{1,(I,J\!-\!\{j^*\})}$}}
\rput(0,-4){\smsize{$\ov\cM_{1,I\sqcup J}$}}
\rput(9,-4){\smsize{$\ov\cM_{1,I\sqcup(J-\{j^*\})}$}}
\psline[linestyle=dashed]{->}(2.2,0)(6.2,0)
\psline{->}(2,-4)(6,-4)\psline{->}(0,-.8)(0,-3.2)\psline{->}(9,-.8)(9,-3.2)
\rput(4.5,.6){\smsize{$\ti{f}$}}\rput(4.5,-3.4){\smsize{$f$}}
\rput(-.5,-2){\smsize{$\pi$}}\rput(8.5,-2){\smsize{$\pi$}}
\end{pspicture}
\caption{Illustration of Main Statement of Proposition~\ref{lift_prp}}
\label{lift_fig1}
\end{figure}

\begin{lmm}
\label{psirestr_lmm1}
With notation as in Proposition~\ref{lift_prp}, for all $i\!\in\!I$
\begin{gather*}
\wt\cM_{1,\vr_{ij^*}} = 
\wt\cM_{1,((I-\{i\})\sqcup\{p\},J-\{j^*\})}\times\ov\cM_{0,\{q,i,j^*\}}
\qquad\hbox{and}\\
\pi_P\!\circ\!\pi=\pi\!\circ\!\pi_P\!:
\wt\cM_{1,\vr_{ij^*}}\lra \ov\cM_{1,((I-\{i\})\sqcup\{p\})\sqcup(J-\{j^*\})},
\end{gather*}
where 
$$\pi_P\!:\wt\cM_{1,\vr_{ij^*}} \lra \wt\cM_{1,((I-\{i\})\sqcup\{p\},J-\{j^*\})}$$
is again the projection onto the first component.
Furthermore, if $\ti\psi$ denotes the universal $\psi$-class and $\ti{f}$ is as in 
Proposition~\ref{lift_prp}, then
$$\ti\psi\big|_{\wt\cM_{1,\vr_{ij^*}}} =0
\qquad\hbox{and}\qquad
\big(\ti{f}^*\ti\psi\big)\big|_{\wt\cM_{1,\vr_{ij^*}}} = \pi_P^*\ti\psi.$$
\end{lmm}

\begin{lmm}
\label{psirestr_lmm2}
With notation as in Proposition~\ref{lift_prp}, for all $j\!\in\!J\!-\!\{j^*\}$
\begin{gather*}
\pi^{-1}\big(\ov\cM_{1,\vr_{jj^*}}\big)
=\wt\cM_{1,\vr_{jj^*}} \approx 
\wt\cM_{1,(I,(J-\{j,j^*\})\sqcup\{p\})}\times\ov\cM_{0,\{q,j,j^*\}}
\qquad\hbox{and}\\
\pi_P\!\circ\!\pi=\pi\!\circ\!\pi_P\!:
\wt\cM_{1,\vr_{jj^*}}\lra \ov\cM_{1,I\sqcup((J-\{j,j^*\})\sqcup\{p\})},
\end{gather*}
where 
$$\pi_P\!:\wt\cM_{1,\vr_{jj^*}} \lra \wt\cM_{1,(I,(J-\{j,j^*\})\sqcup\{p\})}$$
is again the projection onto the first component.
Furthermore, if $\ti\psi$ denotes the universal $\psi$-class on $\wt\cM_{1,(I,J)}$
and on $\wt\cM_{1,(I,(J-\{j,j^*\})\sqcup\{p\})}$, then
$$\ti\psi\big|_{\wt\cM_{1,\vr_{jj^*}}} 
= \big(\ti{f}^*\ti\psi\big)\big|_{\wt\cM_{1,\vr_{jj^*}}} = \pi_P^*\ti\psi.$$\\
\end{lmm}

\noindent
We are now ready to verify the first identity in Theorem~\ref{main_thm}.
We can assume that $\ti{c}\!\neq\!0$; otherwise, 
it reduces to the standard genus-one string equation.
Note that if $i_1,i_2\!\in\!I$ and $i_1\!\neq\!i_2$, then
\begin{equation}\label{psiprod_e0}
\ov\cM_{1,\vr_{i_1j^*}} \cap \ov\cM_{1,\vr_{i_2j^*}} =\eset
\qquad \Lra \qquad
\wt\cM_{1,\vr_{i_1j^*}} \cap \wt\cM_{1,\vr_{i_2j^*}} =\eset.
\end{equation}
Thus, by Proposition~\ref{lift_prp} and Lemma~\ref{psirestr_lmm1},
applied repeatedly,
\begin{equation}\label{psiprod_e1}
\ti\psi^{\ti{c}}=\ti\psi^{\ti{c}-1}\big(\ti{f}^*{\psi}+
\sum_{i\in I}\wt\cM_{1,\vr_{ij^*}}\big)
=\ti{f}^*\ti\psi^{\ti{c}}+
\sum_{i\in I}\big(\pi_P^*\ti\psi^{\ti{c}-1}\big)\!\cap\!\wt\cM_{1,\vr_{ij^*}}.
\end{equation}
On the other hand, by~\e_ref{psipullback_e1b} and Lemma~\ref{psirestr_lmm2},
\begin{equation}\label{psiprod_e2}
\pi^*\psi_j^{c_j}=\ti{f}^*\pi^*\psi_j^{c_j}
+\big(\pi_P^*\pi^*\psi_p^{c_j-1}\big)\!\cap\!\wt\cM_{1,\vr_{jj^*}}
\qquad\forall\,j\!\in\!J\!-\!\{j^*\}.
\end{equation}
If $c_j\!=\!0$, we define the last term in~\e_ref{psiprod_e2} to be zero.
Similarly to~\e_ref{psiprod_e0},
\begin{equation}\label{psiprod_e0b}
\ov\cM_{1,\vr_{ij^*}} \cap \ov\cM_{1,\vr_{jj^*}} =\eset
\quad \Lra \quad
\wt\cM_{1,\vr_{ij^*}} \cap \wt\cM_{1,\vr_{jj^*}} =\eset
\qquad\forall\,j\!\in\!J\!-\!\{j^*\},\,i\!\in\!I\!\sqcup\!J\!-\!\{j,j^*\}.
\end{equation}
Thus, by~\e_ref{psiprod_e1}, \e_ref{psiprod_e2}, and 
Lemmas~\ref{psirestr_lmm1} and~\ref{psirestr_lmm2},
\begin{equation*}\begin{split}
&\blr{\ti{c};(c_j)_{j\in J-\{j^*\}}}_{(|I|,|J|)} 
\equiv \Big\lan \ti\psi^{\ti{c}}
\cdot\!\!\prod_{j\in J-\{j^*\}}\!\!\!\!\!\!\!\pi^*\psi_j^{c_j},\wt\cM_{1,(I,J)}\Big\ran\\
&\qquad\qquad= \Big\lan \ti{f}^*\Big(\ti\psi^{\ti{c}}
\cdot\!\!\prod_{j\in J-\{j^*\}}\!\!\!\!\!\!\!\pi^*\psi_j^{c_j}\Big),\wt\cM_{1,(I,J)}\Big\ran
+\sum_{i\in I} \Big\lan \pi_P^*\Big(\ti\psi^{\ti{c}-1}
\cdot\!\!\prod_{j\in J-\{j^*\}}\!\!\!\!\!\!\!\pi^*\psi_j^{c_j}\Big),
\wt\cM_{1,\vr_{ij^*}}\Big\ran\\
&\qquad\qquad\qquad\qquad\qquad\qquad\qquad\qquad
+\!\!\sum_{j\in J-\{j^*\}}\!\!\!\!\!\!
\Big\lan \pi_P^*\Big(\ti\psi^{\ti{c}}\cdot\pi^*\psi_p^{c_j-1}
\cdot\!\!\prod_{j'\in J-\{j^*,j\}}\!\!\!\!\!\!\!\!\!\!\pi^*\psi_{j'}^{c_{j'}}\Big),
\wt\cM_{1,\vr_{jj^*}}\Big\ran\\
&\qquad\qquad= 0
+\sum_{i\in I} \Big\lan\ti\psi^{\ti{c}-1}
\cdot\!\!\prod_{j\in J-\{j^*\}}\!\!\!\!\!\!\!\pi^*\psi_j^{c_j},
\wt\cM_{1,((I-\{i\})\sqcup\{p\},J-\{j^*\})}\Big\ran\\
&\qquad\qquad\qquad\qquad\qquad\qquad
+\!\!\sum_{j\in J-\{j^*\}}\!\!\!\!\!\!
\Big\lan \ti\psi^{\ti{c}}\cdot\pi^*\psi_p^{c_j-1}
\cdot\!\!\prod_{j'\in J-\{j^*,j\}}\!\!\!\!\!\!\!\!\!\!\pi^*\psi_{j'}^{c_{j'}},
\wt\cM_{1,(I,(J-\{j,j^*\})\sqcup\{p\})}\Big\ran\\
&\qquad\qquad\equiv |I|\blr{\ti{c}\!-\!1;(c_j)_{j\in J-\{j^*\}}}_{(|I|,|J|-1)}
+\!\sum_{j\in J-\{j^*\}}\!\!\!\!\!\!
      \blr{\ti{c};c_j\!-\!1,(c_{j'})_{j'\in J-\{j^*,j\}}}_{(|I|,|J|-1)}.
\end{split}\end{equation*}
We have thus derived the first identity in Theorem~\ref{main_thm}
from Proposition~\ref{lift_prp} and Lemmas~\ref{psirestr_lmm1} and~\ref{psirestr_lmm2}.

\section{Proofs of Main Structural Results}
\label{pf_sec}

\subsection{Proof of Lemma~\ref{psirestr_lmm1}}
\label{lmmpf_subs1}

\noindent
Suppose $I$ is a finite set and $i,j$ are distinct elements of $I$.
It is well-known that the normal bundle $\N_{\ov\cM_{1,I}}\ov\cM_{1,\vr_{ij}}$
of $\ov\cM_{1,\vr_{ij}}$ in $\ov\cM_{1,I}$ is given~by
\begin{equation}\label{normalbundle_e0}
\N_{\ov\cM_{1,I}}\ov\cM_{1,\vr_{ij}}=\pi_P^*L_p\!\otimes\!\pi_B^*L_q
=\pi_P^*L_p,
\end{equation}
where $\pi_B$ and $\pi_B$ are as in~\e_ref{projmap_e} and
$$L_p\lra\ov\cM_{1,(I-\{i,j\})\sqcup\{p\}} \qquad\hbox{and}\qquad 
L_q\lra\ov\cM_{0,\{q,i,j\}}$$
are the universal tangent line bundles for the marked points~$p$ and~$q$;
see \cite{P}, for example.
The last equality in~\e_ref{normalbundle_e0} is due to the fact that 
$\ov\cM_{0,\{q,i,j\}}$ consists of one point.\\

\noindent
Suppose in addition that 
\begin{equation}\label{vrdfn_e}
\vr\!\equiv\!\big(I_P,\{I_k\!:k\!\in\!K\}\big)\in\A_1(I)
\qquad\hbox{and}\qquad \vr\prec\vr_{ij}.
\end{equation}
Then, by the definition of the partial ordering $\prec$ in~\e_ref{partorder_e},
$$\{i,j\}\subset I_k \qquad\hbox{for some}\qquad k\in K.$$
We define $\mu_{ij}(\vr)\!\in\!\A_1\big((I\!-\!\{i,j\})\!\sqcup\!\{p\}\big)$ by
\begin{equation}\label{mudfn_e}
\mu_{ij}(\vr) =\begin{cases}
\big(I_P\!\sqcup\!\{p\},\{I_{k'}\!:k'\!\in\!K-\!\{k\}\}\big),&
\hbox{if}~I_k\!=\!\{i,j\};\\
\big(I_P,\big\{(I_k\!-\!\{i,j\})\!\sqcup\!\{p\}\big\}\!\sqcup\!
\{I_{k'}\!:k'\!\in\!K-\!\{k\}\}\big),& \hbox{if}~I_k\!\supsetneq\!\{i,j\}.
\end{cases}
\end{equation}
It is straightforward to see that 
\begin{equation}\label{muinter_e}
\ov\cM_{1,\vr_{ij}} \cap \ov\cM_{1,\vr}
=\ov\cM_{1,\mu_{ij}(\vr)}\times\ov\cM_{0,\{q,i,j\}}
\subset\ov\cM_{1,(I-\{i,j\})\sqcup\{p\}} \times \ov\cM_{0,\{q,i,j\}}.
\end{equation}

\begin{lmm}
\label{psirestr_lmm1a} 
If $I$ and $J$ are finite sets, $i\!\in\!I$, and $j\!\in\!J$, then the map
\begin{equation}\label{psirestr_lmm1a_e}
\mu_{ij}\!: \big\{\vr\!\in\!\A_1(I,J)\!: \vr\!\prec\!\vr_{ij}\big\}
\lra \A_1\big((I\!-\!\{i\})\!\sqcup\!\{p\},J\!-\!\{j\}\big)
\end{equation}
is an isomorphism of partially ordered sets.
\end{lmm}

\noindent
This lemma follows easily from \e_ref{partorder_e} and~\e_ref{mudfn_e}.
It implies that given an order $<$ on 
$$\A_1\big((I\!-\!\{i\})\!\sqcup\!\{p\},J\!-\!\{j\}\big)$$
extending the partial ordering~$\prec$, we can choose an order $<$ on $A_1(I,J)$ 
that extends the partial ordering $\prec$ such~that
$$\vr_1,\vr_2\prec\vr_{ij}, \quad \mu_{ij}(\vr_1)<\mu_{ij}(\vr_2)
\qquad\Lra\qquad \vr_1<\vr_2.$$
Below we refer to the constructions of Subsection~\ref{blowup_subs} for 
the sets 
$$\A_1\big((I\!-\!\{i\})\!\sqcup\!\{p\},J\!-\!\{j\}\big)
\qquad\hbox{and}\qquad \A_1(I,J)$$
corresponding to such compatible orders~$<$.
We extend the map $\mu_{ij}$ of~\e_ref{psirestr_lmm1a_e} to $\{0\}\!\sqcup\!\A_1(I,J)$
by setting
$$\mu_{ij}(\vr) = \begin{cases}
\mu_{ij}(\max\{\vr'\!<\!\vr:\vr'\!\prec\!\vr_{ij}\}),&
\hbox{if}~\exists\, \vr'\!<\!\vr~\hbox{s.t.}~\vr'\!\prec\!\vr_{ij};\\
0,&\hbox{otherwise}.\end{cases}$$

\begin{lmm}
\label{psirestr_lmm1b} 
Suppose $I$ and $J$ are finite sets, $i\!\in\!I$, and $j\!\in\!J$.
If $\vr\!\in\!\A_1(I,J)$ and $\vr\!<\!\vr_{ij}$, then with notation 
as in Subsection~\ref{blowup_subs} and in~\e_ref{decomp_e}
\begin{gather*}
\ov\cM_{1,\vr_{ij}}^{\vr} = 
\ov\cM_{1,((I-\{i\})\sqcup\{p\}),J-\{j\})}^{\mu_{ij}(\vr)}\times\ov\cM_{0,\{q,i,j\}},\\
\E_{\vr}\big|_{\ov\cM_{1,\vr_{ij}}^{\vr}} = \pi_P^*\E_{\mu_{ij}(\vr)},
\qquad\hbox{and}\qquad
\N_{\ov\cM_{1,(I,J)}^{\vr}}\ov\cM_{1,\vr_{ij}}^{\vr}
=\pi_P^*L_{\mu_{ij}(\vr),p},
\end{gather*}
where 
$$\pi_P\!: \ov\cM_{1,((I-\{i\})\sqcup\{p\},J-\{j\})}^{\mu_{ij}(\vr)}
\times\ov\cM_{0,\{q,i,j\}} \lra
\ov\cM_{1,((I-\{i\})\sqcup\{p\},J-\{j\})}^{\mu_{ij}(\vr)}$$
is the projection map onto the first component.
\end{lmm}

\noindent
By \e_ref{decomp_e}, \e_ref{larestr_e1}, and \e_ref{normalbundle_e0}, 
Lemma~\ref{psirestr_lmm1b}  holds for $\vr\!=\!0$.
Suppose $\vr\!\in\!\A_1(I,J)$, $\vr\!<\!\vr_{ij}$, and the three claims hold for $\vr\!-\!1$.
If $\vr\!\not\prec\!\vr_{ij}$, then
\begin{gather}
\mu_{ij}(\vr)=\mu_{ij}(\vr\!-\!1) \qquad\Lra\notag\\
\label{prtransform_e1a}
\ov\cM_{1,((I-\{i\})\sqcup\{p\},J-\{j\})}^{\mu_{ij}(\vr)}
=\ov\cM_{1,((I-\{i\})\sqcup\{p\},J-\{j\})}^{\mu_{ij}(\vr-1)}, ~~
\E_{\mu_{ij}(\vr)}=\E_{\mu_{ij}(\vr-1)},
~~
L_{\mu_{ij}(\vr),p}=L_{\mu_{ij}(\vr-1),p}.
\end{gather}
On the other hand, since $\vr$ and $\vr_{ij}$ are not comparable with respect to~$\prec$,
the blowup locus $\ov\cM_{1,\vr}^{\vr-1}$ in $\ov\cM_{1,(I,J)}^{\vr-1}$
is disjoint from $\ov\cM_{1,\vr_{ij}}^{\vr-1}$; see Subsection~\ref{blowup_subs} above
and Lemma~\ref{g1desing-curve1bl_lmm} in~\cite{VZ}.
Thus, 
\begin{equation}\label{prtransform_e1b}
\ov\cM_{1,\vr_{ij}}^{\vr}=\ov\cM_{1,\vr_{ij}}^{\vr-1}, \quad
\E_{\vr}|_{\ov\cM_{1,\vr_{ij}}^{\vr}}=\E_{\vr-1}|_{\ov\cM_{1,\vr_{ij}}^{\vr-1}},
\quad 
\N_{\ov\cM_{1,(I,J)}^{\vr}}\ov\cM_{1,\vr_{ij}}^{\vr}
=\N_{\ov\cM_{1,(I,J)}^{\vr-1}}\ov\cM_{1,\vr_{ij}}^{\vr-1}.
\end{equation}
By \e_ref{prtransform_e1a}, \e_ref{prtransform_e1b}, and the inductive assumptions,
the three claims hold for~$\vr$.\\

\noindent
Suppose that $\vr\!\prec\!\vr_{ij}$.
Since all varieties $\ov\cM_{1,\vr'}$ intersect properly in $\ov\cM_{1,(I,J)}$
in the sense of Subsection~\ref{g1desing-curveprelim_subs} in~\cite{VZ},
so do their proper transforms $\ov\cM_{1,\vr'}^{\vr-1}$ in $\ov\cM_{1,(I,J)}^{\vr-1}$.
Furthermore, 
$$\ov\cM_{1,\vr_{ij}}^{\vr-1}\!\cap\!\ov\cM_{1,\vr}^{\vr-1}
\subset \ov\cM_{1,\vr_{ij}}^{\vr-1}\subset \ov\cM_{1,(I,J)}^{\vr-1}$$
is the proper transform of 
$$\ov\cM_{1,\vr_{ij}}\!\cap\!\ov\cM_{1,\vr}
\subset \ov\cM_{1,\vr_{ij}} \subset \ov\cM_{1,(I,J)}.$$
Since $\vr\!\prec\!\vr_{ij}$, $\mu_{ij}(\vr\!-\!1)\!=\!\mu_{ij}(\vr)\!-\!1$.
Thus, by~\e_ref{muinter_e} and the inductive assumptions, 
$$\ov\cM_{1,\vr_{ij}}^{\vr-1}\!\cap\!\ov\cM_{1,\vr}^{\vr-1}
=\ov\cM_{1,\mu_{ij}(\vr)}^{\mu_{ij}(\vr)-1}
\!\times\!\ov\cM_{0,\{q,i,j\}}
\subset \ov\cM_{1,((I-\{i\})\sqcup\{p\},J-\{j\})}^{\mu_{ij}(\vr)-1}
\!\times\!\ov\cM_{0,\{q,i,j\}}.$$
Since $\ov\cM_{1,\vr_{ij}}^{\vr-1}$ and $\ov\cM_{1,\vr}^{\vr-1}$ intersect properly,
the proper transform of $\ov\cM_{1,\vr_{ij}}^{\vr-1}$ in $\ov\cM_{1,(I,J)}^{\vr}$,
i.e.~the blowup of $\ov\cM_{1,(I,J)}^{\vr-1}$ along $\ov\cM_{1,\vr}^{\vr-1}$,
is the blowup of $\ov\cM_{1,\vr_{ij}}^{\vr-1}$ along
$\ov\cM_{1,\vr_{ij}}^{\vr-1}\!\cap\!\ov\cM_{1,\vr}^{\vr-1}$; 
see Subsection~\ref{g1desing-curveprelim_subs} in~\cite{VZ}.
Thus, $\ov\cM_{1,\vr_{ij}}^{\vr}$ is the blowup of 
$$\ov\cM_{1,((I-\{i\})\sqcup\{p\},J-\{j\})}^{\mu_{ij}(\vr)-1}\times\ov\cM_{0,\{q,i,j\}}$$
along $\ov\cM_{1,\mu_{ij}(\vr)}^{\mu_{ij}(\vr)-1}\!\times\!\ov\cM_{0,\{q,i,j\}}$.
By the construction of Subsection~\ref{blowup_subs}, this blowup~is
$$\ov\cM_{1,((I-\{i\})\sqcup\{p\},J-\{j\})}^{\mu_{ij}(\vr)}\times\ov\cM_{0,\{q,i,j\}}.$$
Furthermore, by~\e_ref{bundtwist_e} and the inductive assumptions,
\begin{equation*}\begin{split}
\E_{\vr}|_{\ov\cM_{1,\vr_{ij}}^{\vr}}
&=\big(\ti\pi_{\vr}^*\E_{\vr-1}\!+\!\ov\cM_{1,\vr}^{\vr}\big)\big|_{\ov\cM_{1,\vr_{ij}}^{\vr}}
=\ti\pi_{\vr}^*\pi_P^*\E_{\mu_{ij}(\vr)-1}
+\ov\cM_{1,\mu_{ij}(\vr)}^{\mu_{ij}(\vr)}\!\times\!\ov\cM_{0,\{q,i,j\}}\\
&=\pi_P^*\big(\ti\pi_{\mu_{ij}(\vr)}^*\E_{\mu_{ij}(\vr)-1}
+\ov\cM_{1,\mu_{ij}(\vr)}^{\mu_{ij}(\vr)}\big)
=\E_{\mu_{ij}(\vr)}.
\end{split}\end{equation*}
We have thus verified two of the three inductive assumptions.\\

\noindent
It remains to determine the normal bundle 
$\N_{\ov\cM_{1,(I,J)}^{\vr}}\ov\cM_{1,\vr_{ij}}^{\vr}$
of $\ov\cM_{1,\vr_{ij}}^{\vr}$ in $\ov\cM_{1,(I,J)}^{\vr}$.
We note that by~\e_ref{bundtwist_e} and~\e_ref{mudfn_e},
\begin{equation}\label{bundtwist_e2}
L_{\mu_{ij}(\vr),p}
=\begin{cases}
\ti\pi_{\mu_{ij}(\vr)-1}^*L_{\mu_{ij}(\vr)-1,p}
\!\otimes\!\O(-\ov\cM_{1,\mu_{ij}(\vr)}^{\mu_{ij}(\vr)}),&
\hbox{if}~I_k\!=\!\{i,j\};\\
\ti\pi_{\mu_{ij}(\vr)-1}^*L_{\mu_{ij}(\vr)-1,p}
,& \hbox{if}~I_k\!\supsetneq\!\{i,j\},
\end{cases}
\end{equation}
if $\vr$ is as in~\e_ref{vrdfn_e}.
Furthermore, if $I_k\!=\!\{i,j\}$, then 
$$\ov\cM_{1,\vr}\subset \ov\cM_{1,\vr_{ij}} \qquad\Lra\qquad
\ov\cM_{1,\vr}^{\vr-1}\subset \ov\cM_{1,\vr_{ij}}^{\vr-1}.$$
Thus, by Subsection~\ref{g1desing-map0prelim_subs1} in~\cite{VZ},
\begin{equation}\label{normalbundle_e1a}\begin{split}
\N_{\ov\cM_{1,(I,J)}^{\vr}}\ov\cM_{1,\vr_{ij}}^{\vr}
&=\ti\pi_{\vr}^*~\N_{\ov\cM_{1,(I,J)}^{\vr-1}}\ov\cM_{1,\vr_{ij}}^{\vr-1}
\otimes\O\big(-\ov\cM_{1,\vr_{ij}}^{\vr}\!\cap\!\ov\cM_{1,\vr}^{\vr}\big)\\
&=\ti\pi_{\vr}^*~\N_{\ov\cM_{1,(I,J)}^{\vr-1}}\ov\cM_{1,\vr_{ij}}^{\vr-1}
\otimes\pi_P^*\O\big(-\ov\cM_{1,\mu_{ij}(\vr)}^{\mu_{ij}(\vr)}\big)
\qquad\hbox{if}\quad I_k\!=\!\{i,j\}.
\end{split}\end{equation}
On the other hand, if $I_k\!\supsetneq\!\{i,j\}$, 
$\ov\cM_{1,\vr}^{\vr-1}$ and  $\ov\cM_{1,\vr_{ij}}^{\vr-1}$ intersect transversally in
$\ov\cM_{1,(I,J)}^{\vr-1}$, since 
$\ov\cM_{1,\vr}$ and  $\ov\cM_{1,\vr_{ij}}$ intersect transversally in $\ov\cM_{1,(I,J)}$.
Thus,
\begin{equation}\label{normalbundle_e1b}
\N_{\ov\cM_{1,(I,J)}^{\vr}}\ov\cM_{1,\vr_{ij}}^{\vr}
=\ti\pi_{\vr}^*~\N_{\ov\cM_{1,(I,J)}^{\vr-1}}\ov\cM_{1,\vr_{ij}}^{\vr-1}
\qquad\hbox{if}\quad I_k\!\supsetneq\!\{i,j\}.
\end{equation}
The final inductive assumption now follows from the corresponding statement for $\vr\!-\!1$,
along with \e_ref{bundtwist_e2}-\e_ref{normalbundle_e1b}.

\begin{crl}
\label{psirestr_crl1}
With notation as in Lemma~\ref{psirestr_lmm1},
\begin{gather*}
\wt\cM_{1,\vr_{ij^*}} = 
\wt\cM_{1,((I-\{i\})\sqcup\{p\},J-\{j^*\})}\times\ov\cM_{0,\{q,i,j^*\}},\\
\ti\psi\big|_{\wt\cM_{1,\vr_{ij^*}}} =0,
\qquad\hbox{and}\qquad
\big(\ti{f}^*\ti\psi\big)\big|_{\wt\cM_{1,\vr_{ij^*}}} = \pi_P^*\ti\psi.
\end{gather*}\\
\end{crl}

\noindent
By the paragraph preceding Proposition~\ref{lift_prp} and the first statement of
Lemma~\ref{psirestr_lmm1b}
\begin{equation*}\begin{split}
\wt\cM_{1,\vr_{ij^*}} = \ov\cM_{1,\vr_{ij^*}}^{\vr_{ij^*}-1} 
&= \ov\cM_{1,((I-\{i\})\sqcup\{p\},J-\{j^*\})}^{\mu_{ij^*}(\vr_{ij^*}-1)}
\times\ov\cM_{0,\{q,i,j^*\}}\\
&=\wt\cM_{1,((I-\{i\})\sqcup\{p\},J-\{j^*\})}\times\ov\cM_{0,\{q,i,j^*\}},
\end{split}\end{equation*}
since  $\mu_{ij^*}(\vr_{ij^*}\!-\!1)$ is the largest element of 
$$\big(\A_1((I\!-\!\{i\})\!\sqcup\!\{p\},J\!-\!\{j^*\}),<\big),$$
according to Lemma~\ref{psirestr_lmm1a}.\\

\noindent
Since $\vr_{ij^*}$ is a maximal element of $(\A_1(I,J),\prec)$,
$$\ov\cM_{1,\vr_{ij^*}}^{\vr-1}\cap \ov\cM_{1,\vr}^{\vr-1} 
=\eset \qquad\forall \vr\in\A_1(I,J),\,\vr>\vr_{ij^*}.$$
Thus, by \e_ref{bundtwist_e} and the second statement of Lemma~\ref{psirestr_lmm1b},
\begin{equation}\label{larest_e5a}
\ti\E\big|_{\wt\cM_{1,\vr_{ij^*}}} 
=\E_{\vr_{ij^*}-1}\big|_{\wt\cM_{1,\vr_{ij^*}}}
+\sum_{\vr\ge\vr_{ij^*}}\!\!\!\ov\cM_{1,\vr}^{\vr}\big|_{\wt\cM_{1,\vr_{ij^*}}}
=\pi_P^*\ti\psi+\wt\cM_{1,\vr_{ij^*}}\big|_{\wt\cM_{1,\vr_{ij^*}}}.
\end{equation}
By the third statement of Lemma~\ref{psirestr_lmm1b},
\begin{equation}\label{larest_e5b}\begin{split}
\wt\cM_{1,\vr_{ij^*}}\big|_{\wt\cM_{1,\vr_{ij^*}}}
=\N_{\wt\cM_{1,(I,J)}}\wt\cM_{1,\vr_{ij^*}}
&=\N_{\ov\cM_{1,(I,J)}^{\vr_{ij^*}-1}}\ov\cM_{1,\vr_{ij^*}}^{\vr_{ij^*}-1}\\
&=\pi_P^*L_{\mu_{ij^*}(\vr_{ij^*}-1),p}
=-\pi_P^*\ti\psi.
\end{split}\end{equation}
The second statement of Corollary~\ref{psirestr_crl1} follows from 
\e_ref{larest_e5a} and~\e_ref{larest_e5b}.\\

\noindent
Finally, by the last statement of Proposition~\ref{lift_prp}, the second statement
of Corollary~\ref{psirestr_crl1}, \e_ref{psiprod_e0}, and~\e_ref{larest_e5b},
$$\big(\ti{f}^*\ti\psi\big)\big|_{\wt\cM_{1,\vr_{ij^*}}} =
\ti\psi|_{\wt\cM_{1,\vr_{ij^*}}} -
\sum_{i'\in I}\wt\cM_{1,\vr_{i'j^*}}\big|_{\wt\cM_{1,\vr_{ij^*}}}
=0-\wt\cM_{1,\vr_{ij^*}}\big|_{\wt\cM_{1,\vr_{ij^*}}}
=\pi_P^*\ti\psi.$$
This concludes the proof of Corollary~\ref{psirestr_crl1}.

\subsection{Proof of Lemma~\ref{psirestr_lmm2}}
\label{lmmpf_subs2}

\noindent
The proof of Lemma~\ref{psirestr_lmm2} is analogous to the previous subsection.
If $I$ is a finite set and $j,j^*$ are distinct elements of~$I$, let
\begin{equation*}\begin{split}
\A_1(I;jj^*)
&=\big\{\vr\!\in\!\A_1(I)-\!\{\vr_{jj^*}\}\!: 
\ov\cM_{1,\vr_{jj^*}}\!\cap\!\ov\cM_{1,\vr}\neq\eset\big\}\\
&=\big\{\big(I_P,\{I_k\!:k\!\in\!K\}\big)\!\in\!\A_1(I)-\!\{\vr_{jj^*}\}\!\!:
\{j,j^*\}\!\subset\!I_k~\hbox{for some}~k\!\in\!\{P\}\!\sqcup\!K\big\}.
\end{split}\end{equation*}
For each $\vr\!\in\!\A_1(I;jj^*)$, we define 
$\eta_{jj^*}(\vr)\!\in\!\A_1\big((I\!-\!\{j,j^*\})\!\sqcup\!\{p\}\big)$ by
\begin{equation}\label{etadfn_e}
\eta_{jj^*}(\vr) =\begin{cases}
\big((I_P\!-\!\{j,j^*\})\!\sqcup\!\{p\},\{I_{k'}\!:k'\!\in\!K\}\big),&
\hbox{if}~I_P\!=\!\{j,j^*\};\\
\big(I_P\!\sqcup\!\{p\},\{I_{k'}\!:k'\!\in\!K-\!\{k\}\}\big),&
\hbox{if}~I_k\!=\!\{j,j^*\};\\
\big(I_P,\big\{(I_k\!-\!\{j,j^*\})\!\sqcup\!\{p\}\big\}\!\sqcup\!
\{I_{k'}\!:k'\!\in\!K-\!\{k\}\}\big),& \hbox{if}~I_k\!\supsetneq\!\{j,j^*\}.
\end{cases}
\end{equation}
It is straightforward to see that 
\begin{equation}\label{etainter_e}
\ov\cM_{1,\vr_{jj^*}} \cap \ov\cM_{1,\vr}
=\ov\cM_{1,\eta_{jj^*}(\vr)}\times\ov\cM_{0,\{q,j,j^*\}}
\subset\ov\cM_{1,(I-\{j,j^*\})\sqcup\{p\}} \times \ov\cM_{0,\{q,j,j^*\}}.
\end{equation}

\begin{lmm}
\label{psirestr_lmm2a} 
If $I$ and $J$ are finite sets, $j,j^*\!\in\!J$, and $j\!\neq\!j^*$, then the map
\begin{equation}\label{psirestr_lmm2a_e}
\eta_{jj^*}\!: \A_1(I,J)\!\cap\!\A_1(I\!\sqcup\!J;jj^*) \lra 
\A_1\big((I,(J\!-\!\{j,j^*\})\!\sqcup\!\{p\}\big)
\end{equation}
is an isomorphism of partially ordered sets.
\end{lmm}

\noindent
This lemma follows easily from~\e_ref{partorder_e} and~\e_ref{etadfn_e}.
Note, however, that it is essential that $j,j^*\!\in\!J$ and thus
the second case of~\e_ref{etadfn_e} does not occur if 
$$\vr \in\!\A_1(I,J)\!\cap\!\A_1(I\!\sqcup\!J;jj^*).$$
Lemma~\ref{psirestr_lmm2a} implies that given an order $<$ on 
$$\A_1\big((I,(J\!-\!\{j,j^*\})\!\sqcup\!\{p\}\big)$$
extending the partial ordering~$\prec$, we can choose an order $<$ on $A_1(I,J)$ 
that extends the partial ordering $\prec$ such~that
$$\vr_1,\vr_2\in\A_1(I,J)\!\cap\!\A_1(I\!\sqcup\!J;jj^*), \quad 
\eta_{jj^*}(\vr_1)<\eta_{jj^*}(\vr_2)
\qquad\Lra\qquad \vr_1<\vr_2.$$
Below we refer to the constructions of Subsection~\ref{blowup_subs} for 
the sets 
$$\A_1\big((I,(J\!-\!\{j,j^*\})\!\sqcup\!\{p\}\big)
\qquad\hbox{and}\qquad \A_1(I,J)$$
corresponding to such compatible orders~$<$.
We extend the map $\eta_{jj^*}$ of~\e_ref{psirestr_lmm2a_e} to $\{0\}\!\sqcup\!\A_1(I,J)$
by setting
$$\eta_{jj^*}(\vr) = \begin{cases}
\eta_{jj^*}(\max\{\vr'\!<\!\vr:\vr'\!\in\!\A_1(I\!\sqcup\!J;jj^*)\}),&
\hbox{if}~\exists\, \vr'\!<\!\vr~\hbox{s.t.}~\vr'\!\in\!\A_1(I\!\sqcup\!J;jj^*);\\
0,&\hbox{otherwise}.\end{cases}$$

\begin{lmm}
\label{psirestr_lmm2b} 
Suppose $I$ and $J$ are finite sets, $j,j^*\!\in\!J$, and $j\!\neq\!j^*$.
If $\vr\!\in\!\A_1(I,J)$, then with notation 
as in Subsection~\ref{blowup_subs} and in~\e_ref{decomp_e}
$$\pi_{\vr}^{-1}\big(\ov\cM_{1,\vr_{jj^*}}\big)=\ov\cM_{1,\vr_{jj^*}}^{\vr} = 
\ov\cM_{1,((I,(J-\{j,j^*\})\sqcup\{p\})}^{\eta_{jj^*}(\vr)}\times\ov\cM_{0,\{q,j,j^*\}},
\qquad
\E_{\vr}\big|_{\ov\cM_{1,\vr_{jj^*}}^{\vr}} = \pi_P^*\E_{\eta_{jj^*}(\vr)},$$
where
$$\pi_P: \ov\cM_{1,((I,(J-\{j,j^*\})\sqcup\{p\})}^{\eta_{jj^*}(\vr)}
\times\ov\cM_{0,\{q,j,j^*\}} \lra
\ov\cM_{1,((I,(J-\{j,j^*\})\sqcup\{p\})}^{\eta_{jj^*}(\vr)}$$
is the projection map onto the first component.
\end{lmm}

\noindent
By \e_ref{decomp_e} and \e_ref{larestr_e1}, 
Lemma~\ref{psirestr_lmm2b}  holds for $\vr\!=\!0$.
Suppose $\vr\!\in\!\A_1(I,J)$ and the three claims hold for $\vr\!-\!1$.
If $\vr\!\not\in\!\A_1(I\!\sqcup\!J,jj^*)$, then
\begin{gather}
\eta_{jj^*}(\vr)=\eta_{jj^*}(\vr\!-\!1) \qquad\Lra\notag\\
\label{prtransform_e2a}
\ov\cM_{1,((I,(J-\{j,j^*\})\sqcup\{p\})}^{\eta_{jj^*}(\vr)}
=\ov\cM_{1,((I,(J-\{j,j^*\})\sqcup\{p\})}^{\eta_{jj^*}(\vr-1)}, \qquad
\E_{\eta_{jj^*}(\vr)}=\E_{\eta_{jj^*}(\vr-1)}.
\end{gather}
On the other hand, since 
$$\ov\cM_{1,\vr_{jj^*}}\!\cap\!\ov\cM_{1,\vr}=\eset,$$
the blowup locus $\ov\cM_{1,\vr}^{\vr-1}$ in $\ov\cM_{1,(I,J)}^{\vr-1}$
is disjoint from $\ov\cM_{1,\vr_{jj^*}}^{\vr-1}$.
Thus, 
\begin{equation}\label{prtransform_e2b}
\pi_{\vr}^{-1}\big(\ov\cM_{1,\vr_{jj^*}}\big)=\pi_{\vr-1}^{-1}\big(\ov\cM_{1,\vr_{jj^*}}\big),
\quad \ov\cM_{1,\vr_{jj^*}}^{\vr}=\ov\cM_{1,\vr_{jj^*}}^{\vr-1}, 
\quad \E_{\vr}|_{\ov\cM_{1,\vr_{jj^*}}^{\vr}}=\E_{\vr-1}|_{\ov\cM_{1,\vr_{jj^*}}^{\vr-1}}.
\end{equation}
By \e_ref{prtransform_e2a}, \e_ref{prtransform_e2b}, and the inductive assumptions,
the three claims hold for~$\vr$.\\

\noindent
Suppose that $\vr\!\in\!\A_1(I\!\sqcup\!J,jj^*)$.
Since all varieties $\ov\cM_{1,\vr'}$ intersect properly in $\ov\cM_{1,(I,J)}$,
so do their proper transforms $\ov\cM_{1,\vr'}^{\vr-1}$, 
with $\vr'\!>\!\vr\!-\!1$, in $\ov\cM_{1,(I,J)}^{\vr-1}$.
Since $\ov\cM_{1,\vr}^{\vr-1}$ is not contained in the divisor $\ov\cM_{1,\vr_{jj^*}}^{\vr-1}$,
$\ov\cM_{1,\vr}^{\vr-1}$ and $\ov\cM_{1,\vr_{jj^*}}^{\vr-1}$ intersect transversally.
Thus, using the first statement of the lemma with $\vr$ replaced by $\vr\!-\!1$, we obtain
$$\pi_{\vr}^{-1}\big(\ov\cM_{1,\vr_{jj^*}}\big)
=\ti\pi_{\vr}^{-1}\pi_{\vr-1}^{-1}\big(\ov\cM_{1,\vr_{jj^*}}\big)
=\ti\pi_{\vr}^{-1}\big(\ov\cM_{1,\vr_{jj^*}}^{\vr-1}\big)
=\ov\cM_{1,\vr_{jj^*}}^{\vr}.$$
Furthermore, 
$$\ov\cM_{1,\vr_{jj^*}}^{\vr-1}\!\cap\!\ov\cM_{1,\vr}^{\vr-1}
\subset \ov\cM_{1,\vr_{jj^*}}^{\vr-1}\subset \ov\cM_{1,(I,J)}^{\vr-1}$$
is the proper transform of 
$$\ov\cM_{1,\vr_{jj^*}}\!\cap\!\ov\cM_{1,\vr}
\subset \ov\cM_{1,\vr_{jj^*}} \subset \ov\cM_{1,(I,J)}.$$
Since $\vr\!\in\!\A_1(I\!\sqcup\!J,jj^*)$, $\eta_{jj^*}(\vr\!-\!1)\!=\!\eta_{jj^*}(\vr)\!-\!1$.
Thus, by~\e_ref{etainter_e} and the inductive assumptions, 
$$\ov\cM_{1,\vr_{jj^*}}^{\vr-1}\!\cap\!\ov\cM_{1,\vr}^{\vr-1}
=\ov\cM_{1,\eta_{jj^*}(\vr)}^{\eta_{jj^*}(\vr)-1}
\!\times\!\ov\cM_{0,\{q,j,j^*\}}
\subset \ov\cM_{1,(I,(J-\{j,j^*\})\sqcup\{p\})}^{\eta_{jj^*}(\vr)-1}
\!\times\!\ov\cM_{0,\{q,j,j^*\}}.$$
Since $\ov\cM_{1,\vr_{jj^*}}^{\vr-1}$ and $\ov\cM_{1,\vr}^{\vr-1}$ intersect properly,
the proper transform of $\ov\cM_{1,\vr_{jj^*}}^{\vr-1}$ in $\ov\cM_{1,(I,J)}^{\vr}$,
i.e.~the blowup of $\ov\cM_{1,(I,J)}^{\vr-1}$ along $\ov\cM_{1,\vr}^{\vr-1}$,
is the blowup of $\ov\cM_{1,\vr_{jj^*}}^{\vr-1}$ along
$\ov\cM_{1,\vr_{jj^*}}^{\vr-1}\!\cap\!\ov\cM_{1,\vr}^{\vr-1}$; 
see Subsection~\ref{g1desing-curveprelim_subs} in~\cite{VZ}.
Thus, $\ov\cM_{1,\vr_{jj^*}}^{\vr}$ is the blowup of 
$$\ov\cM_{1,(I,(J-\{j,j^*\})\sqcup\{p\})}^{\eta_{jj^*}(\vr)-1}\times\ov\cM_{0,\{q,j,j^*\}}$$
along $\ov\cM_{1,\eta_{jj^*}(\vr)}^{\eta_{jj^*}(\vr)-1}\!\times\!\ov\cM_{0,\{q,j,j^*\}}$.
By the construction of Subsection~\ref{blowup_subs}, this blowup~is
$$\ov\cM_{1,(I,(J-\{j,j^*\})\sqcup\{p\})}^{\eta_{jj^*}(\vr)}\times\ov\cM_{0,\{q,j,j^*\}}.$$
Furthermore, by~\e_ref{bundtwist_e} and the inductive assumptions,
\begin{equation*}\begin{split}
\E_{\vr}|_{\ov\cM_{1,\vr_{jj^*}}^{\vr}}
&=\big(\ti\pi_{\vr}^*\E_{\vr-1}\!+\!\ov\cM_{1,\vr}^{\vr}\big)\big|_{\ov\cM_{1,\vr_{jj^*}}^{\vr}}
=\ti\pi_{\vr}^*\pi_P^*\E_{\eta_{jj^*}(\vr)-1}
+\ov\cM_{1,\eta_{jj^*}(\vr)}^{\eta_{jj^*}(\vr)}\!\times\!\ov\cM_{0,\{q,j,j^*\}}\\
&=\pi_P^*\big(\ti\pi_{\eta_{jj^*}(\vr)}^*\E_{\eta_{jj^*}(\vr)-1}
+\ov\cM_{1,\eta_{jj^*}(\vr)}^{\eta_{jj^*}(\vr)}\big)
=\E_{\eta_{jj^*}(\vr)}.
\end{split}\end{equation*}
We have thus verified the three inductive assumptions.

\begin{crl}
\label{psirestr_crl2}
With notation as in Proposition~\ref{lift_prp},
\begin{gather*}
\pi^{-1}\big(\ov\cM_{1,\vr_{jj^*}}\big)
=\wt\cM_{1,\vr_{jj^*}} \approx 
\wt\cM_{1,(I,(J-\{j,j^*\})\sqcup\{p\})}\times\ov\cM_{0,\{q,j,j^*\}}\\
\hbox{and}\qquad
\ti\psi\big|_{\wt\cM_{1,\vr_{jj^*}}} 
= \big(\ti{f}^*\ti\psi\big)\big|_{\wt\cM_{1,\vr_{jj^*}}} = \pi_P^*\ti\psi.
\end{gather*}\\
\end{crl}

\noindent
By Lemma~\ref{psirestr_lmm2a}, $\eta_{jj^*}(\vr_{\max})$ is the largest element of 
$$\big(\A_1(I,(J\!-\!\{j,j^*\})\!\sqcup\!\{p\}),<\big).$$
Thus, by the first two statements of Lemma~\ref{psirestr_lmm2b},
\begin{equation*}\begin{split}
\pi^{-1}\big(\ov\cM_{1,\vr_{jj^*}}\big)
&=\pi_{\vr_{\max}}^{-1}\big(\ov\cM_{1,\vr_{jj^*}}\big)
=\ov\cM_{1,\vr_{jj^*}}^{\vr_{\max}}
=\wt\cM_{1,\vr_{jj^*}}\\
&= \ov\cM_{1,(I,(J-\{j,j^*\})\sqcup\{p\})}^{\eta_{jj^*}(\vr_{\max})}
\times\ov\cM_{0,\{q,j,j^*\}}
=\wt\cM_{1,(I,(J-\{j,j^*\})\sqcup\{p\})}\times\ov\cM_{0,\{q,j,j^*\}}.
\end{split}\end{equation*}
By the last statement of Lemma~\ref{psirestr_lmm2b},
$$\ti\psi\big|_{\wt\cM_{1,\vr_{jj^*}}}
=\E_{\vr_{\max}}\big|_{\wt\cM_{1,\vr_{jj^*}}}
=\pi_P^*\E_{\eta_{jj^*}(\vr_{\max})}
=\pi_P^*\ti\E=\pi_P^*\ti\psi.$$
Finally, by the last statement of Proposition~\ref{lift_prp} and \e_ref{psiprod_e0b},
$$\big(\ti{f}^*\ti\psi\big)\big|_{\wt\cM_{1,\vr_{jj^*}}} =
\ti\psi|_{\wt\cM_{1,\vr_{jj^*}}} -
\sum_{i\in I}\wt\cM_{1,\vr_{ij^*}}\big|_{\wt\cM_{1,\vr_{jj^*}}}
=\pi_P^*\ti\psi-0.$$
This concludes the proof of Corollary~\ref{psirestr_crl2}.

\subsection{Proof of Proposition~\ref{lift_prp}}
\label{lift_subs}

\noindent
In this subsection we prove Proposition~\ref{lift_prp}.
In fact, we show that there is a lift of the forgetful map~$f$
of Proposition~\ref{lift_prp} to morphisms between corresponding
stages of the blowup construction of Subsection~\ref{blowup_subs}
for $\ov\cM_{1,(I,J)}$ and for $\ov\cM_{1,(I,J-\{j^*\})}$;
see Lemma~\ref{lift_lmm} below.\\

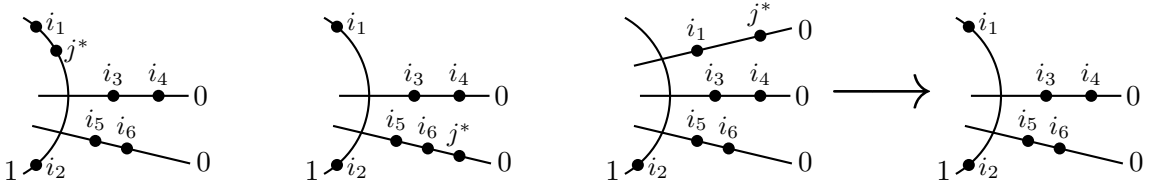
\begin{figure}
\begin{pspicture}(-1.1,-1.8)(10,1.3)
\psset{unit=.4cm}
\psarc(-2,-1){3}{-60}{60}\pscircle*(-.07,1.3){.2}\pscircle*(-.07,-3.3){.2}
\rput(.6,1.4){\begin{small}$i_1$\end{small}}
\rput(.6,-3.4){\begin{small}$i_2$\end{small}}
\psline(0,-1)(5,-1)\pscircle*(2.5,-1){.2}\pscircle*(4,-1){.2}
\pscircle*(.6,.5){.2}\rput(1.3,.5){\begin{small}$j^*$\end{small}}
\rput(2.5,-.3){\begin{small}$i_3$\end{small}}
\rput(4,-.3){\begin{small}$i_4$\end{small}}
\psline(-.2,-2)(5.05,-3.25)\pscircle*(1.9,-2.5){.2}
\pscircle*(2.95,-2.75){.2}
\rput(1.9,-1.8){\begin{small}$i_5$\end{small}}
\rput(2.95,-2.05){\begin{small}$i_6$\end{small}}
\rput(-.9,-3.5){$1$}\rput(5.4,-1){$0$}\rput(5.5,-3.2){$0$}
\psarc(8,-1){3}{-60}{60}\pscircle*(9.93,1.3){.2}\pscircle*(9.93,-3.3){.2}
\rput(10.6,1.4){\begin{small}$i_1$\end{small}}
\rput(10.6,-3.4){\begin{small}$i_2$\end{small}}
\psline(10,-1)(15,-1)\pscircle*(12.5,-1){.2}\pscircle*(14,-1){.2}
\rput(12.5,-.3){\begin{small}$i_3$\end{small}}
\rput(14,-.3){\begin{small}$i_4$\end{small}}
\psline(9.8,-2)(15.05,-3.25)\pscircle*(11.9,-2.5){.2}
\pscircle*(12.95,-2.75){.2}\pscircle*(14,-3){.2}
\rput(11.9,-1.8){\begin{small}$i_5$\end{small}}
\rput(12.95,-2.05){\begin{small}$i_6$\end{small}}
\rput(14,-2.3){\begin{small}$j^*$\end{small}}
\rput(9.1,-3.5){$1$}\rput(15.4,-1){$0$}\rput(15.5,-3.2){$0$}
\psarc(18,-1){3}{-60}{60}\pscircle*(19.93,-3.3){.2}
\rput(20.6,-3.4){\begin{small}$i_2$\end{small}}
\psline(19.8,0)(25.05,1.25)\pscircle*(21.9,.5){.2}\pscircle*(24,1){.2}
\rput(21.9,1.2){\begin{small}$i_1$\end{small}}
\rput(24,1.7){\begin{small}$j^*$\end{small}}\psline(20,-1)(25,-1)\pscircle*(22.5,-1){.2}\pscircle*(24,-1){.2}
\rput(22.5,-.3){\begin{small}$i_3$\end{small}}
\rput(24,-.3){\begin{small}$i_4$\end{small}}
\psline(19.8,-2)(25.05,-3.25)\pscircle*(21.9,-2.5){.2}
\pscircle*(22.95,-2.75){.2}
\rput(21.9,-1.8){\begin{small}$i_5$\end{small}}
\rput(22.95,-2.05){\begin{small}$i_6$\end{small}}
\rput(19.1,-3.5){$1$}\rput(25.5,1.2){$0$}\rput(25.4,-1){$0$}\rput(25.5,-3.2){$0$}
\rput(28,-1){\begin{Huge}$\lra$\end{Huge}}
\psarc(29,-1){3}{-60}{60}\pscircle*(30.93,1.3){.2}\pscircle*(30.93,-3.3){.2}
\rput(31.6,1.4){\begin{small}$i_1$\end{small}}
\rput(31.6,-3.4){\begin{small}$i_2$\end{small}}
\psline(31,-1)(36,-1)\pscircle*(33.5,-1){.2}\pscircle*(35,-1){.2}
\rput(33.5,-.3){\begin{small}$i_3$\end{small}}
\rput(35,-.3){\begin{small}$i_4$\end{small}}
\psline(30.8,-2)(36.05,-3.25)\pscircle*(32.9,-2.5){.2}
\pscircle*(33.95,-2.75){.2}
\rput(32.9,-1.8){\begin{small}$i_5$\end{small}}
\rput(33.95,-2.05){\begin{small}$i_6$\end{small}}
\rput(30.1,-3.5){$1$}\rput(36.4,-1){$0$}\rput(36.5,-3.2){$0$}
\end{pspicture}
\caption{Images under the Forgetful Map}
\label{forgetfulmap_fig}
\end{figure}

\noindent
First, we define a forgetful map 
$$f\!:\A_1(I,J)\lra \bar\A_1\big(I,J\!-\!\{j^*\}\big)\equiv 
\A_1\big(I,J\!-\!\{j^*\}\big) \sqcup \big\{(I\!\sqcup\!(J\!-\!\{j^*\}),\eset)\big\}.$$ 
If $\vr\!=\!(I_P\!\sqcup\!J_P,\{I_k\!\sqcup\!J_k\!: k\!\in\!K\})$, we put
$$f(\vr) = \begin{cases}
\big(I_P\!\sqcup\!(J_P\!-\!\{j^*\}),\{I_k\!\sqcup\!J_k\!: k\!\in\!K\}\big),&
\hbox{if}~j^*\!\in\!J_P;\\
\big(I_P\!\sqcup\!J_P,\{I_k\!\sqcup\!(J_k\!-\!\{j^*\})\}\!\sqcup\!
\{I_{k'}\!\sqcup\!J_{k'}\!: k'\!\in\!K\!-\!\{k\}\}\big),
&\hbox{if}~j^*\!\in\!J_k,~|I_k|\!+\!|J_k|\!>\!2;\\
\big((I_P\!\sqcup\!\{i\})\!\sqcup\!J_P,\{I_{k'}\!\sqcup\!J_{k'}\!: k'\!\in\!K\!-\!\{k\}\}\big),
&\hbox{if}~I_k\!\sqcup\!J_k\!=\!\{ij^*\}.
\end{cases}$$
These three cases are represented in Figure~\ref{forgetfulmap_fig}.
We note that for all $\rho\!\in\!\A_1(I,J\!-\!\{j^*\})$,
$$f^{-1}\big(\ov\cM_{1,\rho}\big)=
\!\!\bigcup_{\vr\in f^{-1}(\rho)}\!\!\!\!\!\ov\cM_{1,\vr}.$$
Furthermore, 
$$\rho_1,\rho_2\!\in\!\bar\A_1(I,J\!-\!\{j^*\}), ~~\rho_1\!\neq\!\rho_2, ~~~
\vr_1\!\in\!f^{-1}(\rho_1), ~~~ \vr_2\!\in\!f^{-1}(\rho_2), ~~~
\vr_1\!\prec\!\vr_2 \quad\Lra\quad \rho_1\!\prec\!\rho_2.$$
Thus, given an order $<$ on $\A_1(I,J\!-\!\{j^*\})$ extending the partial ordering $\prec$,
we can choose an order $<$ on $\A_1(I,J)$ extending $\prec$ such~that
$$\rho_1,\rho_2\!\in\!\bar\A_1(I,J\!-\!\{j^*\}), ~~\rho_1\!<\!\rho_2, ~~~
\vr_1\!\in\! f^{-1}(\rho_1), ~~~ \vr_2\!\in\!f^{-1}(\rho_2)
 \quad\Lra\quad \vr_1\!<\!\vr_2.$$
Below we will refer to the blowup constructions of Subsection~\ref{blowup_subs}
for $\ov\cM_{1,(I,J)}$ and for $\ov\cM_{1,(I,J-\{j^*\})}$
corresponding to such compatible orders.
For each $\rho\!\in\!\A_1(I,J\!-\!\{j^*\})$, let
$$\rho^+=\max f^{-1}(\rho)\in\A_1(I,J) \qquad\hbox{and}\qquad
\rho^-=\min f^{-1}(\rho)-1\in\{0\}\!\sqcup\!\A_1(I,J).$$
If $\rho$ is not the minimal element of $\A_1(I,J\!-\!\{j^*\})$,
then $\rho^-\!=\!(\rho\!-\!1)^+$.

\begin{lmm}
\label{lift_lmm}
Suppose $I$, $J$, and $f$ are as in Proposition~\ref{lift_prp}.
For each $\rho\!\in\!\A_1(I,J\!-\!\{j^*\})$, $f$ lifts to a morphism
$$f_{\rho}\!:\ov\cM_{1,(I,J)}^{\rho^+} \lra \ov\cM_{1,(I,J-\{j^*\})}^{\rho}$$
over the projection maps 
$$\pi_{\rho^+}\!: \ov\cM_{1,(I,J)}^{\rho^+} \lra \ov\cM_{1,(I,J)}
\qquad\hbox{and}\qquad
\pi_{\rho}\!: \ov\cM_{1,(I,J-\{j^*\})}^{\rho} \lra \ov\cM_{1,(I,J-\{j^*\})};$$
see Figure~\ref{lift_fig2}.
Furthermore, 
\begin{equation}\label{lift_lmm_e}
f_{\rho}^{-1}\big(\ov\cM_{1,\rho^*}^{\rho}\big)
=\!\!\bigcup_{\vr\in f^{-1}(\rho^*)}\!\!\!\!\!\!\ov\cM_{1,\vr}^{\rho^+}
\quad\forall\,\rho^*\!>\!\rho \qquad\hbox{and}\qquad
\E_{\rho^+}=f_{\rho}^*\E_{\rho}.
\end{equation}\\
\end{lmm}

\begin{figure}
\begin{pspicture}(-2,-1.8)(10,1.3)
\psset{unit=.4cm}
\rput(0,0){\smsize{$\ov\cM_{1,(I,J)}^{\rho^+}$}}
\rput(9,0){\smsize{$\ov\cM_{1,(I,J\!-\!\{j^*\})}^{\rho}$}}
\rput(0,-4){\smsize{$\ov\cM_{1,I\sqcup J}$}}
\rput(9,-4){\smsize{$\ov\cM_{1,I\sqcup(J-\{j^*\})}$}}
\psline[linestyle=dashed]{->}(2.2,0)(6.2,0)
\psline{->}(2,-4)(6,-4)\psline{->}(0,-.8)(0,-3.2)\psline{->}(9,-.8)(9,-3.2)
\rput(4.5,.6){\smsize{$f_{\rho}$}}\rput(4.5,-3.4){\smsize{$f$}}
\rput(-.8,-2){\smsize{$\pi_{\rho^+}$}}\rput(8.4,-2){\smsize{$\pi_{\rho}$}}
\rput(20,0){\smsize{$\ov\cM_{1,(I,J)}^{\rho^+}$}}
\rput(29,0){\smsize{$\ov\cM_{1,(I,J\!-\!\{j^*\})}^{\rho}$}}
\rput(20,-4){\smsize{$\ov\cM_{1,(I,J)}^{\rho^-}$}}
\rput(29,-4){\smsize{$\ov\cM_{1,(I,J\!-\!\{j^*\})}^{\rho-1}$}}
\psline[linestyle=dashed]{->}(22.2,0)(26.2,0)
\psline{->}(22,-4)(26,-4)\psline{->}(20,-.8)(20,-3.2)\psline{->}(29,-.8)(29,-3.2)
\rput(24.5,.6){\smsize{$f_{\rho}$}}\rput(24.5,-3.4){\smsize{$f_{\rho-1}$}}
\rput(16.7,-2){\smsize{$\ti\pi_{\rho^-+1}\circ\ldots\circ\ti\pi_{\rho^+}$}}
\rput(28.4,-2){\smsize{$\ti\pi_{\rho}$}}
\end{pspicture}
\caption{Main Statement of Lemma~\ref{lift_lmm} and Inductive Step in the Proof}
\label{lift_fig2}
\end{figure}

\noindent
Proposition~\ref{lift_prp} follows easily from Lemma~\ref{lift_lmm}
by taking $\rho\!=\!\rho_{\max}$,
where $\rho_{\max}$ is the maximal element of $\A_1(I,J\!-\!\{j^*\})$.
We note that
$$\big\{\vr\!\in\!\A_1(I,J)\!: \vr\!>\!\rho_{\max}^+\big\}
=\big\{\vr\!\in\!\A_1(I,J)\!: f(\vr)\!=\!(I\!\sqcup\!(J\!-\!\{j^*\}),\eset)\big\}
=\big\{\vr_{ij^*}\!: i\!\in\!I\big\}.$$
Since $\ov\cM_{1,\vr_{ij^*}}\!\subset\!\ov\cM_{1,I\sqcup J}$ is a divisor 
for every $i\!\in\!I$, so is 
$$\ov\cM_{1,\vr_{ij^*}}^{\rho_{\max}^+} \subset\ov\cM_{1,I\sqcup J}^{\rho_{\max}^+}.$$
Thus, by the construction of Subsection~\ref{blowup_subs},
\begin{gather*}
\wt\cM_{1,I\sqcup J} \equiv \ov\cM_{1,I\sqcup J}^{\vr_{\max}}
=\ov\cM_{1,I\sqcup J}^{\rho_{\max}^+} \qquad\hbox{and}\\
\E\equiv\E_{\vr_{\max}}=\E_{\rho_{\max}^+} +\sum_{i\in I}\ov\cM_{1,\vr_{ij^*}}^{\rho_{\max}^+} 
=f_{\rho_{\max}}^*\E_{\rho_{\max}}+\sum_{i\in I}\ov\cM_{1,\vr_{ij^*}}^{\rho_{\max}^+}
=\ti{f}^*\E+\sum_{i\in I}\wt\cM_{1,\vr_{ij^*}},
\end{gather*}
where $\ti{f}\!=\!f_{\rho_{\max}}$.\\

\noindent
Lemma~\ref{lift_lmm} holds for $\rho\!=\!0\!\in\!\{0\}\!\cup\!\A_1(I,J\!-\!\{j^*\})$,
if we define $0^+\!=\!0$.
Suppose 
$$\rho\!=\!\big(I_P\!\sqcup\!J_P,\{I_k\!\sqcup\!J_k\!:k\!\in\!K\}\big)   
\in\A_1(I,J\!-\!\{j^*\})$$ 
and the lemma holds~for 
$$\rho\!-\!1\in\{0\}\!\sqcup\!\A_1(I,J\!-\!\{j^*\}).$$
The elements of $f^{-1}(\rho)\!\subset\!\A_1(I,J)$ can be described as follows.
The largest element~is
$$\rho^+=\big(I_P\!\sqcup\!(J_P\!\sqcup\!\{j^*\}),\{I_k\!\sqcup\!J_k\!:k\!\in\!K\}\big).$$
Furthermore, for each $k\!\in\!K$ and $i\!\in\!I_P$,
\begin{alignat*}{1}
\rho_k(j^*) &\equiv \big(I_P\!\sqcup\!J_P,
\big\{I_k\!\sqcup\!(J_k\!\sqcup\!\{j^*\})\big\}\!\sqcup\!
\{I_{k'}\!\sqcup\!J_{k'}\!:k'\!\in\!K\!-\!\{k\}\}\big)\in f^{-1}(\rho);\\
\rho_i(j^*) &\equiv \big((I_P\!-\!\{i\})\!\sqcup\!J_P,
\big\{\{i,j\}\big\}\!\sqcup\!
\{I_{k'}\!\sqcup\!J_{k'}\!:k'\!\in\!K\}\big)\in f^{-1}(\rho).
\end{alignat*}
It is straightforward to see that 
$$f^{-1}(\rho)=\big\{\rho_k(j^*)\!: k\!\in\!K\big\} \sqcup
\big\{\rho_i(j^*)\!: i\!\in\!I_P\big\} \sqcup 
 \{\rho^+\}.$$
Furthermore, $\rho^+$ is the largest element of $(f^{-1}(\rho),\prec)$,
while no two elements of the form $\rho_k(j^*)$ and/or $\rho_i(j^*)$
are comparable with respect to~$\prec$. Thus, 
$$\ov\cM_{1,\rho_k(j^*)}^{\rho^-} \cap  \ov\cM_{1,\rho_i(j^*)}^{\rho^-} =\eset
\qquad\forall\, i,k\in I_P\!\sqcup\!K,\,i\!\neq\!k;$$
see Subsection~\ref{blowup_subs}.
In fact,
$$\ov\cM_{1,\rho_k(j^*)} \cap  \ov\cM_{1,\rho_i(j^*)} =\eset
\qquad\forall\, i,k\in I_P\!\sqcup\!K,\,i\!\neq\!k;$$
see the proof of Lemma~\ref{g1desing-curve1bl_lmm} in~\cite{VZ}.\\

\noindent
All varieties $\ov\cM_{1,\vr^*}$ are smooth and intersect properly in $\ov\cM_{1,(I,J)}$
in the sense of Subsection~\ref{g1desing-curveprelim_subs} in~\cite{VZ}.
Thus, all varieties $\ov\cM_{1,\vr^*}^{\rho^-}$, with $\vr^*\!>\!\rho^-$,
are also smooth and intersect properly in $\ov\cM_{1,(I,J)}^{\rho^-}$.
It follows that for every $k\!\in\!K$ and every point 
$$p\in  \ov\cM_{1,\rho_k(j^*)}^{\rho^-}-  \ov\cM_{1,\rho^+}^{\rho^-},$$ 
we can choose neighborhoods $\ti{U}$ of $p$ in $\ov\cM_{1,(I,J)}^{\rho^-}$,
$U$ of $f_{\rho-1}(p)$ in $\ov\cM_{1,(I,J-\{j^*\})}^{\rho-1}$, and  
coordinates~$(z,v,t)$ on~$\ti{U}$ such~that\\
${}\quad~~$ (i) $U\!=\!f_{\rho-1}(\ti{U})$;\\
${}\quad~$ (ii) $U\!=\!\big\{(z,v)\!\in\!\C^{|I|+|J|-|K|-1}\!\times\!\C^K\big\}$;\\
${}\quad$ (iii) $\ov\cM_{1,\rho}^{\rho-1}\!\cap\!U\!=\!
\big\{(z,v)\!\in\!U\!: v\!=\!0\big\}$;\\
${}\quad$ (iv) $\ti{U}\!=\!\big\{(z,v,t)\!\in\! 
\C^{|I|+|J|-|K|-1}\!\times\!\C^K\!\times\!\C\big\}$
and $f_{\rho-1}(z,v,t)=(z,v)$.\\
These assumptions imply that 
$$\ov\cM_{1,\rho_k(j^*)}^{\rho^-}\!\cap\!\ti{U}=
\big\{(z,v,t)\!\in\!\ti{U}\!: v\!=\!0\big\}.$$
Since $\ov\cM_{1,(I,J-\{j^*\})}^{\rho}$ is the blowup of $\ov\cM_{1,(I,J-\{j^*\})}^{\rho-1}$
along $\ov\cM_{1,\rho}^{\rho-1}$, the preimage of~$U$ in $\ov\cM_{1,(I,J-\{j^*\})}^{\rho}$
under the projection map~is
$$V=\big\{(z,v;\ell)\!\in\!U\!\times\!\bP\big(\C^K\big)\!:v\!\in\!\ell\big\}.$$
Since $\ov\cM_{1,(I,J)}^{\rho^+}$ is the blowup of $\ov\cM_{1,(I,J)}^{\rho^-}$
along $\ov\cM_{1,\rho_k(j^*)}^{\rho^-}$ and subvarieties that do not contain~$p$,
the preimage of~$\ti{U}$ in $\ov\cM_{1,(I,J)}^{\rho^+}$
under the projection map~is
$$\ti{V}=\big\{(z,v,t;\ell)\!\in\!\ti{U}\!\times\!\bP\big(\C^K\big)\!:v\!\in\!\ell\big\},$$
provided $\ti{U}$ is sufficiently small.
Thus, the map $f_{\rho-1}\!:\ti{U}\!\lra\!U$ lifts to a~map $f_{\rho}\!:\ti{V}\!\lra\!V$.
This lift is defined~by
\begin{equation}\label{lift_e1}
 f_{\rho}(z,v,t;\ell) = (z,v;\ell).
\end{equation}\\

\noindent
Similarly to the previous paragraph, for every 
$$p\in   \ov\cM_{1,\rho^+}^{\rho^-}-
 \bigcup_{k\in K}\!\ov\cM_{1,\rho_k(j^*)}^{\rho^-},$$ 
we can choose neighborhoods $\ti{U}$ of $p$ in $\ov\cM_{1,(I,J)}^{\rho^-}$,
$U$ of $f_{\rho-1}(p)$ in $\ov\cM_{1,(I,J-\{j^*\})}^{\rho-1}$, and  
coordinates~$(z,v,t)$ on~$\ti{U}$ such~that the conditions (i)-(iv) are satisfied,
with $\ov\cM_{1,\rho_k(j^*)}^{\rho^-}$ replaced by $\ov\cM_{1,\rho^+}^{\rho^-}$.
Thus, if 
$$p \not\in \bigcup_{i\in I_P}\!\ov\cM_{1,\rho_i(j^*)}^{\rho^-},$$
the map $f_{\rho-1}$ lifts to the preimage of a neighborhood of $p$ in 
$\ov\cM_{1,(I,J)}^{\rho^+}$, just as in the previous paragraph.\\

\noindent
On the other hand, suppose 
$$p \in\ov\cM_{1,\rho_i(j^*)}^{\rho^-}$$
for some $i\!\in\!I_P$. 
Since $\ov\cM_{1,\rho_i(j^*)}\!\subset\!\ov\cM_{1,\rho^+}$
is of codimension-one, 
$$\ov\cM_{1,\rho_i(j^*)}^{\rho^-}\subset\ov\cM_{1,\rho^+}^{\rho^-}$$
is also of codimension-one.
We can thus choose local coordinate so that
$$\ov\cM_{1,\rho_i(j^*)}^{\rho^-}\cap\ti{U}
=\big\{(z,v,t)\!\in\!\ti{U}\!: v\!=\!0,\, t\!=\!0\big\}.$$
Since $\ov\cM_{1,(I,J)}^{\rho^+-1}$ is the blowup of $\ov\cM_{1,(I,J)}^{\rho^-}$
along $\ov\cM_{1,\rho_i(j^*)}^{\rho^-}$ and subvarieties that do not contain~$p$,
the preimage of~$\ti{U}$ in $\ov\cM_{1,(I,J)}^{\rho^+-1}$
under the projection map~is
$$\ti{V}=\big\{(z,v,t;\ell')\!\in\!\ti{U}\!\times\!\bP\big(\C^K\!\times\!\C\big): 
(v,t)\!\in\!\ell'\big\},$$
provided $\ti{U}$ is sufficiently small.
It is immediate that
$$ \ov\cM_{1,\rho^+}^{\rho^+-1}\cap \ti{V}
=\big\{(z,0,t;[\al,\be])\!\in\!\ti{U}\!\times\!\bP\big(\C^K\!\times\!\C\big)\!: 
\al\!=\!0\big\},$$
where $\ov\cM_{1,\rho^+}^{\rho^+-1}\!\subset\!\ov\cM_{1,(I,J)}^{\rho^+-1}$ 
is the proper transform of $\ov\cM_{1,\rho^+}^{\rho^-}$.
A neighborhood of $\ov\cM_{1,\rho^+}^{\rho^+-1}\!\cap\!\ti{V}$ is given~by
$$\ti{U}'=\big\{(z,u,t)\!\in\C^{|I|+|J|-|K|-1}\!\times\!\C^K\!\times\!\C\big\},
\qquad (z,u,t) \llra \big(z,ut,t;[u,1]\big)\in\ti{V}.$$
Since $\ov\cM_{1,(I,J)}^{\rho^+}$ is the blowup of $\ov\cM_{1,(I,J)}^{\rho^+-1}$
along $\ov\cM_{1,\rho^+}^{\rho^+-1}$,
the preimage of~$\ti{U}$ in $\ov\cM_{1,(I,J)}^{\rho^+}$
under the projection map~is 
\begin{gather*}
\wt{W}=\big(\big\{(z,u,t;\ell)\!\in\!\ti{V}'\!\times\!\bP\big(\C^K\big)\!:
u\!\in\!\ell\big\}
\cup\big\{(z,v,t;[\al,\be])\!\in\!\ti{V}\!:\al\!\neq\!0\big\}\big)\big/\sim,\\
(z,u,t;\ell) \sim  \big(z,ut,t;[u,1]\big).
\end{gather*}
Thus, the map $f_{\rho-1}\!:\ti{U}\!\lra\!U$ lifts to a~map $f_{\rho}\!:\wt{W}\!\lra\!V$.
This lift is defined~by
\begin{equation}\label{lift_e2}
f_{\rho}(z,u,t;\ell) = (z,ut;\ell) \qquad\hbox{and}\qquad 
f_{\rho}\big(z,v,t;[\al,\be]\big)=\big(z,v;[\al]\big)
\end{equation}
on the two charts on~$\wt{W}$.
Note that if $u\!\neq\!0$, then $[u]\!=\!\ell\!\in\!\bP(\C^K)$.
Thus, the map~$f_{\rho}$ agrees on the overlap of the two charts.\\

\noindent
Finally, suppose that 
$$p\in \ov\cM_{1,\rho_k(j^*)}^{\rho^-}\cap \ov\cM_{1,\rho^+}^{\rho^-}$$
for some $k\!\in\!K$.
Since the varieties $\ov\cM_{1,\vr^*}$ intersect properly in $\ov\cM_{1,(I,J)}$,
$\ov\cM_{1,\rho_k(j^*)}^{\rho^-}$ and  $\ov\cM_{1,\rho^+}^{\rho^-}$ intersect 
properly in~$\ov\cM_{1,(I,J)}^{\rho^-}$ and
$\ov\cM_{1,\rho_k(j^*)}^{\rho^-}\!\cap\!\ov\cM_{1,\rho^+}^{\rho^-}$ is
the proper transform of $\ov\cM_{1,\rho_k(j^*)}\!\cap\!\ov\cM_{1,\rho^+}$.
Thus, $\ov\cM_{1,\rho_k(j^*)}^{\rho^-}\!\cap\!\ov\cM_{1,\rho^+}^{\rho^-}$
is a divisor in  $\ov\cM_{1,\rho_k(j^*)}^{\rho^-}$ and in $\ov\cM_{1,\rho^+}^{\rho^-}$.
If follows that we can choose neighborhoods $\ti{U}$ of $p$ in $\ov\cM_{1,(I,J)}^{\rho^-}$,
$U$ of $f_{\rho-1}(p)$ in $\ov\cM_{1,(I,J-\{j^*\})}^{\rho-1}$, and  
coordinates $(z,v,w_k,w_+)$ on~$\ti{U}$ such~that\\
${}\quad~~$ (i) $U\!=\!f_{\rho-1}(\ti{U})$;\\
${}\quad~$ (ii) $U\!=\!\big\{(z,v,w)\!\in\!\C^{|I|+|J|-|K|-1}
\!\times\!\C^{K-\{k\}}\!\times\!\C\big\}$;\\
${}\quad$ (iii) $\ov\cM_{1,\rho}^{\rho-1}\!\cap\!U\!=\!
\big\{(z,v,w)\!\in\!U\!: v\!=\!0,\,w\!=\!0\big\}$;\\
${}\quad$ (iv) $\ti{U}\!=\!\big\{(z,v,w_k,w_+)\!\in\! \C^{|I|+|J|-|K|-1}\times\!\C^{K-\{k\}}\!\times\!\C\!\times\!\C\big\}$,
$f_{\rho-1}(z,v,w_k,w_+)=(z,v,w_kw_+)$;\\
${}\quad~$ (v) $\ov\cM_{1,\rho_k(j^*)}^{\rho^-}\!\cap\!\ti{U}\!=\!
\big\{(z,v,w_k,w_+)\!\in\!\ti{U}\!: \,v\!=\!0,\,w_+\!=\!0\big\}$;\\
${}\quad$ (vi) $\ov\cM_{1,\rho^+}^{\rho^-}\!\cap\!\ti{U}\!=\!
\big\{(z,v,w_k,w_+)\!\in\!\ti{U}\!: v\!=\!0,\, w_k\!=\!0\big\}$.\\
Similarly to the above, the preimage of~$U$ in $\ov\cM_{1,(I,J-\{j^*\})}^{\rho}$
under the projection map~is
$$V=\big\{(z,v,w;\ell)\!\in\!U\!\times\!\bP\big(\C^{K-\{k\}}\!\times\!\C\big)\!:
(v,w)\!\in\!\ell\big\}.$$
Since $\ov\cM_{1,(I,J)}^{\rho^+-1}$ is the blowup of $\ov\cM_{1,(I,J)}^{\rho^-}$
along $\ov\cM_{1,\rho_k(j^*)}^{\rho^-}$ and subvarieties that do not contain~$p$,
the preimage of~$\ti{U}$ in $\ov\cM_{1,(I,J)}^{\rho^+-1}$
under the projection map~is
$$\ti{V}=\big\{(z,v,w_k,w_+;\ell_k)\!\in\!
\ti{U}\!\times\!\bP\big(\C^{K-\{k\}}\!\times\!\C\big)\!: (v,w_+)\!\in\!\ell_k\big\},$$
provided $\ti{U}$ is sufficiently small.
It is immediate that
$$\ov\cM_{1,\rho^+}^{\rho^+-1} \cap \ti{V}
=\big\{(z,0,0,w_+;[\al,\be])\!\in\!\ti{U}\!\times\bP\big(\C^{K-\{k\}}\!\times\!\C\big)\!: 
\al\!=\!0\big\},$$
where $\ov\cM_{1,\rho^+}^{\rho^+-1}\!\subset\!\ov\cM_{1,(I,J)}^{\rho^+-1}$ 
is the proper transform of $\ov\cM_{1,\rho^+}^{\rho^-}$.
A neighborhood of $\ov\cM_{1,\rho^+}^{\rho^+-1}\!\cap\!\ti{V}$ is given~by
\begin{gather*}
\ti{U}'=\big\{(z,u,u_k,w_+)\!\in\!
\C^{|I|+|J|-|K|-1}\!\times\!\C^{K-\{k\}}\!\times\!\C\!\times\!\C\big\},\\
(z,u,u_k,w_+) \llra \big(z,uw_+,u_k,w_+;[u,1]\big)\in\ti{V}.
\end{gather*}
Since $\ov\cM_{1,(I,J)}^{\rho^+}$ is the blowup of $\ov\cM_{1,(I,J)}^{\rho^+-1}$
along $\ov\cM_{1,\rho^+}^{\rho^+-1}$,
the preimage of~$\ti{U}$ in $\ov\cM_{1,(I,J)}^{\rho^+}$
under the projection map~is 
\begin{gather*}\begin{split}
\wt{W}=&\big(\big\{(z,u,u_k,w_+;\ell)\!\in\!
\ti{V}'\!\times\!\bP\big(\C^{K-\{k\}}\!\times\!\C\big)\!: (u,u_k)\!\in\!\ell\big\}\\
&\qquad\qquad\qquad\qquad
\cup\big\{(z,v,w_k,w_+;[\al,\be])\!\in\!\ti{V}\!:\al\!\neq\!0\big\}\big)\big/\sim,
\end{split}\\
 (z,u,u_k,w_+;\ell) \sim \big(z,uw_+,u_k,w_+;[u,1]\big) .
\end{gather*}
Thus, the map $f_{\rho-1}\!:\ti{U}\!\lra\!U$ lifts to a~map $f_{\rho}\!:\wt{W}\!\lra\!V$.
This lift is defined~by
\begin{equation}\label{lift_e3}\begin{split}
f_{\rho}(z,u,u_k,w_+;\ell) &= (z,uw_+,u_kw_+;\ell) \qquad\hbox{and}\\
f_{\rho}\big(z,v,w_k,w_+;[\al,\be]\big) &=(z,v,w_kw_+;[\al,w_k\be])
\end{split}\end{equation}
on the two charts on~$\wt{W}$. It is immediate that $f_{\rho}$ is well-defined
on the overlap of the two charts.\\

\noindent
{\it Remark:} The first equality in \e_ref{lift_lmm_e} should be viewed
as incorporating the above information concerning the local structure 
of the projection map.
It is easy to see from the verification of the first equality in \e_ref{lift_lmm_e} below
that this additional information is preserved by the inductive step as~well.\\

\noindent
It remains to verify that the two equalities in \e_ref{lift_lmm_e} still hold.
Let 
\begin{gather*}
\pi_{\rho,\rho-1}\!: \ov\cM_{1,(I,J-\{j^*\})}^{\rho} \lra 
 \ov\cM_{1,(I,J-\{j^*\})}^{\rho-1}
\qquad\hbox{and}\\
\pi_{\rho^+,\vr}\!: \ov\cM_{1,(I,J)}^{\rho^+}\lra\ov\cM_{1,(I,J)}^{\vr},
\quad\vr\in \{\rho^-\}\!\cup\!f^{-1}(\rho)
\end{gather*}
be the projection maps.
By the construction of the line bundles $\E_{\vr}$ in Subsection~\ref{blowup_subs},
\begin{gather}\label{lift_lmm_e3a}
\E_{\rho}=\pi_{\rho,\rho-1}^{~*}\E_{\rho}+\ov\cM_{1,\rho}^{\rho}
\qquad\hbox{and}\\
\label{lift_lmm_e3b}
\E_{\rho^+}=\pi_{\rho^+,\rho^-}^{~*}\E_{\rho^-}+
\sum_{\vr\in f^{-1}(\rho)}\!\!\!\!\pi_{\rho^+,\vr}^*\ov\cM_{1,\vr}^{\vr}
=\pi_{\rho^+,\rho^-}^{~*}\E_{\rho^-}+
\sum_{\vr\in f^{-1}(\rho)}\!\!\!\!\pi_{\rho^+,\vr}^{-1}\big(\ov\cM_{1,\vr}^{\vr}\big),
\end{gather}
where
$$\ov\cM_{1,\rho}^{\rho}=\pi_{\rho,\rho-1}^{~-1}\big(\ov\cM_{1,\rho}^{\rho-1}\big)
\subset\ov\cM_{1,(I,J-\{j^*\})}^{\rho} 
\qquad\hbox{and}\qquad
\ov\cM_{1,\vr}^{\vr}\subset \pi_{\vr,\vr-1}^{~-1}\big(\ov\cM_{1,\vr}^{\vr-1}\big)$$
are the exceptional divisors for the blowups at the steps $\rho$ and $\vr$.
Since all divisors $\pi_{\rho^+,\vr}^{-1}\big(\ov\cM_{1,\vr}^{\vr}\big)$
are distinct,
\begin{equation}\label{lift_lmm_e3c}\begin{split}
\sum_{\vr\in f^{-1}(\rho)}\!\!\!\!\pi_{\rho^+,\vr}^{-1}\big(\ov\cM_{1,\vr}^{\vr}\big)
&=\pi_{\rho^+,\rho^-}^{~-1}
 \bigg(\bigcup_{\vr\in f^{-1}(\rho)}\!\!\!\!\!\!\ov\cM_{1,\vr}^{\rho^-}\bigg)
=\pi_{\rho^+,\rho^-}^{~-1}\big(f_{\rho-1}^{\,-1}(\ov\cM_{1,\rho}^{\rho-1})\big)\\
&=f_{\rho}^{\,-1}\pi_{\rho,\rho-1}^{~-1}\big(\ov\cM_{1,\rho}^{\rho-1}\big)
=f_{\rho}^{\,-1}\big(\ov\cM_{1,\rho}^{\rho}\big)
=f_{\rho}^{\,*}\big(\ov\cM_{1,\rho}^{\rho}\big).
\end{split}\end{equation}
The second equality in \e_ref{lift_lmm_e} follows from the same equality with
$\rho$ replaced by $\rho\!-\!1$, along with \e_ref{lift_lmm_e3a}-\e_ref{lift_lmm_e3c}.\\

\noindent
Suppose next that $\rho^*\!>\!\rho$.
Since
$$\pi_{\rho,\rho-1}\circ f_{\rho}=f_{\rho-1}\circ\pi_{\rho^+,\rho^-},$$
$\ov\cM_{1,\rho^*}^{\rho}$ is the proper transform of $\ov\cM_{1,\rho^*}^{\rho-1}$,
and $\ov\cM_{1,\vr^*}^{\rho^+}$ is the proper transform of $\ov\cM_{1,\vr^*}^{\rho^-}$,
$$f_{\rho}^{-1}\big(\ov\cM_{1,\rho^*}^{\rho}\big)
\supset \!\bigcup_{\vr^*\in f^{-1}(\rho^*)}\!\!\!\!\!\!\!\!\ov\cM_{1,\vr^*}^{\rho^+}$$
by the first equation in~\e_ref{lift_lmm_e} with $\rho$ replaced by $\rho\!-\!1$.
We will next verify the opposite inclusion.
Suppose 
\begin{gather*}
q\in\ov\cM_{1,\rho^*}^{\rho}, \qquad \ti{p}\in f_{\rho}^{-1}(q),\qquad\hbox{and}\\
p=\pi_{\rho^+,\rho^-}(\ti{p})
\in f_{\rho-1}^{\,-1}\big(\ov\cM_{1,\rho^*}^{\rho-1}\big)
= \!\bigcup_{\vr^*\in f^{-1}(\rho^*)}\!\!\!\!\!\!\!\!\ov\cM_{1,\vr^*}^{\rho^-}
\subset\ov\cM_{1,(I,J)}^{\rho^-}.
\end{gather*}
If $\pi_{\rho,\rho-1}(q)\!\not\in\!\ov\cM_{1,\rho}^{\rho-1}$, then
$$f_{\rho}^{-1}(q)=f_{\rho-1}^{-1}\big(\pi_{\rho,\rho-1}(q)\big)
=p\in 
\bigcup_{\vr^*\in f^{-1}(\rho^*)}\!\!\!\!\!\!\!\!\ov\cM_{1,\vr^*}^{\rho^-}
-\!\bigcup_{\vr\in f^{-1}(\rho)}\!\!\!\!\!\!\ov\cM_{1,\vr}^{\rho^-}
\subset \bigcup_{\vr^*\in f^{-1}(\rho^*)}\!\!\!\!\!\!\!\!\ov\cM_{1,\vr^*}^{\rho^+},$$
as needed.\\

\noindent
Suppose that 
$$\pi_{\rho,\rho-1}(q) \in 
\ov\cM_{(\rho,\rho^*)}^{\rho-1}\equiv
\ov\cM_{1,\rho}^{\rho-1} \cap \ov\cM_{1,\rho^*}^{\rho-1}.$$
First, we consider the case when
$$p\in \ov\cM_{1,\rho_k(j^*)}^{\rho^-}-\ov\cM_{1,\rho^+}^{\rho^-}$$
for some $k\!\in\!K$.
Since $\ov\cM_{1,\rho}^{\rho-1}$ and $\ov\cM_{1,\rho^*}^{\rho-1}$
intersect properly in $\ov\cM_{1,(I,J-\{j^*\})}^{\rho-1}$, we can choose local coordinates
$(z,v,t)$ near $p$ as in the first case considered above such~that for some 
$K_{\rho^*}\!\subset\!K$\\
${}\quad$ (v) $\ov\cM_{1,\rho^*}^{\rho-1}\!\cap\!U\!=\!
\big\{(z,v)\!\in\!U\!: z\!\in\!\ov\cM_{(\rho,\rho^*)}^{\rho-1};\,
 v\!\in\!\C^{K_{\rho^*}}\big\}$.\\
This assumption implies that 
\begin{equation}\label{baseblow_e1}
\ov\cM_{1,\rho}^{\rho}\cap \ov\cM_{1,\rho^*}^{\rho}\cap V
=\big\{(z,0;\ell)\!\in\!V\!\!: z\!\in\!\ov\cM_{(\rho,\rho^*)}^{\rho-1};
~\ell\!\in\!\bP(\C^{K_{\rho^*}})\big\}.
\end{equation}
In addition, by (iv) and the structure of $f_{\rho-1}$,
$$\bigcup_{\vr^*\in f^{-1}(\rho^*)}\!\!\!\!\!\!\!\!\ov\cM_{1,\vr^*}^{\rho^-}
\cap \ti{U}
=f_{\rho-1}^{-1}\big(\ov\cM_{1,\rho^*}^{\rho-1}\big)\cap\ti{U} 
=\big\{(z,v,t)\!\in\!\ti{U}\!:
z\!\in\!\ov\cM_{(\rho,\rho^*)}^{\rho-1};~
v\!\in\!\C^{K_{\rho^*}}\big\}.$$
Since $\ov\cM_{1,\rho_k(j^*)}^{\rho^-}$ and $\ov\cM_{1,\vr^*}^{\rho^-}$ 
intersect properly, it follows that
$$\bigcup_{\vr^*\in f^{-1}(\rho^*)}\!\!\!\!\!\!\!\!\ov\cM_{1,\vr^*}^{\rho^+}
\cap\ov\cM_{1,\rho_k(j^*)}^{\rho^+}\cap \ti{V}
=\big\{(z,0,t;\ell)\!\in\!\ti{V}; z\!\in\!\ov\cM_{(\rho,\rho^*)}^{\rho-1};
~\ell\!\in\!\bP(\C^{K_{\rho^*}})\big\}.$$
Using~\e_ref{lift_e1}, we conclude that
$$\ti{p} \in \big\{f_{\rho}|_{\ti{V}}\big\}^{-1}
\big(\ov\cM_{1,\rho^*}^{\rho}\!\cap\!\ov\cM_{1,\rho}^{\rho}\big)
=\bigcup_{\vr^*\in f^{-1}(\rho^*)}\!\!\!\!\!\!\!\!\ov\cM_{1,\vr^*}^{\rho^+}
\cap\ov\cM_{1,\rho_k(j^*)}^{\rho^+}\cap \ti{V},$$
as needed.\\

\noindent
Suppose next that
$$p\in \ov\cM_{1,\rho^+}^{\rho^-} - \bigcup_{k\in K}\ov\cM_{1,\rho_k(j^*)}^{\rho^-},$$
i.e.~as in the second case considered above.
We can again choose $K_{\rho^*}\!\subset\!K$ so that Condition~(v)
in the previous paragraph is satisfied.
If
$$p \not\in \bigcup_{i\in I_P}\ov\cM_{1,\rho_i(j^*)}^{\rho^-},$$
then the same argument as in the previous paragraph, but 
with replaced $\rho_k(j^*)$ by $\rho^+$, shows that 
$$\ti{p}\in\bigcup_{\vr^*\in f^{-1}(\rho^*)}\!\!\!\!\!\!\!\!\ov\cM_{1,\vr^*}^{\rho^+}.$$
On the other hand, suppose that
$$p\in \ov\cM_{1,\rho_i(j^*)}^{\rho^-}$$
for some $i\!\in\!I_P$. 
Then, with notation as in the construction of the map $f_{\rho}$ in this case,
\begin{gather*}
\bigcup_{\vr^*\in f^{-1}(\rho^*)}\!\!\!\!\!\!\!\!\ov\cM_{1,\vr^*}^{\rho^+-1} \cap\ti{V} 
=\big\{(z,v,t;\ell')\!\in\!\ti{V}\!:
z\!\in\!\ov\cM_{(\rho,\rho^*)}^{\rho-1};\,
\ell'\!\in\!\bP(\C^{K_{\rho^*}}\!\times\!\C)\big\} ~~\Lra\\
\bigcup_{\vr^*\in f^{-1}(\rho^*)}\!\!\!\!\!\!\!\!\ov\cM_{1,\vr^*}^{\rho^+-1}\cap\ti{U}'  
=\big\{(z,u,t)\!\in\!\ti{U}'\!\!: 
z\!\in\!\ov\cM_{(\rho,\rho^*)}^{\rho-1};\,
u\!\in\!\C^{K_{\rho^*}}\big\} \qquad\Lra\\
\begin{split}
\bigcup_{\vr^*\in f^{-1}(\rho^*)}\!\!\!\!\!\!\!\!\ov\cM_{1,\vr^*}^{\rho^+}
\cap  \pi_{\rho^+,\rho^-}^{~-1}\big(\ov\cM_{1,\rho_i(j^*)}^{\rho^-}\big)
&\cap\wt{W} 
=\big\{(z,u,0;\ell)\!\in\!\wt{W}\!: z\!\in\!\ov\cM_{(\rho,\rho^*)}^{\rho-1};\,
\ell\!\in\!\bP(\C^{K_{\rho^*}})\!\subset\!\bP\big(\C^K\big)\big\}\\
&\cup \big\{(z,0,0;\ell')\!\in\!\wt{W}\!: z\!\in\!\ov\cM_{(\rho,\rho^*)}^{\rho-1};\,
\ell'\!\in\!\bP\big(\C^{K_{\rho^*}}\!\times\!\C\big)\!\subset\!\bP\big(\C^K\!\times\!\C\big)\big\}.
\end{split}
\end{gather*}
Using \e_ref{lift_e2} and~\e_ref{baseblow_e1}, we conclude that
$$\ti{p} \in \big\{
f_{\rho}|_{\pi_{\rho^+,\rho^-}^{~-1}(\ov\cM_{1,\rho_i(j^*)}^{\rho^-})\cap\wt{W}}\big\}^{-1}
\big(\ov\cM_{1,\rho^*}^{\rho}\!\cap\!\ov\cM_{1,\rho}^{\rho}\big)
=\bigcup_{\vr^*\in f^{-1}(\rho^*)}\!\!\!\!\!\!\!\!\ov\cM_{1,\vr^*}^{\rho^+}
\cap  \pi_{\rho^+,\rho^-}^{~-1}\big(\ov\cM_{1,\rho_i(j^*)}^{\rho^-}\big)
\cap\wt{W}.$$
Note that the map 
$f_{\rho}|_{\pi_{\rho^+,\rho^-}^{~-1}(\ov\cM_{1,\rho_i(j^*)}^{\rho^-}) \cap \wt{W}}$ 
is a $\bP^1$-fibration,
while the map $f_{\rho}|_{\ti{V}}$ of the previous paragraph is a $\C$-fibration.\\

\noindent
Finally, suppose that 
$$p\in \ov\cM_{1,\rho^+}^{\rho^-}\cap \ov\cM_{1,\rho_k(j^*)}^{\rho^-}$$
for some $k\!\in\!K$.
With notation as in the corresponding case in the construction of 
the map~$f_{\rho}$ and with a good choice of local coordinates, 
we have two cases to consider.
There exists $K_{\rho^*}\!\subset\!K\!-\!\{k\}$ such~that\\
${}\quad$ Case~1: 
$\ov\cM_{1,\rho^*}^{\rho-1}\!\cap\!U\!=\!
\big\{(z,v,w)\!\in\!U\!: z\!\in\!\ov\cM_{(\rho,\rho^*)}^{\rho-1};
\, v\!\in\!\C^{K_{\rho^*}}\big\}$;\\
${}\quad$ Case~2: 
$\ov\cM_{1,\rho^*}^{\rho-1}\!\cap\!U\!=\!
\big\{(z,v,w)\!\in\!U\!: z\!\in\!\ov\cM_{(\rho,\rho^*)}^{\rho-1};
\, v\!\in\!\C^{K_{\rho^*}},\,w\!=\!0\big\}$.\\
In the first case, we have
\begin{gather}\label{baseblow_e2}
\ov\cM_{1,\rho^*}^{\rho}\cap\ov\cM_{1,\rho}^{\rho}\cap V 
=\big\{(z,0,0;\ell)\!\in\!V\!\!: 
z\!\in\!\ov\cM_{(\rho,\rho^*)}^{\rho-1};
~\ell\!\in\!\bP(\C^{K_{\rho^*}}\!\times\!\C)\big\} \qquad\hbox{and}\\
\bigcup_{\vr^*\in f^{-1}(\rho^*)}\!\!\!\!\!\!\!\!\ov\cM_{1,\vr^*}^{\rho^-} \cap\ti{U} 
=f_{\rho-1}^{-1}\big(\ov\cM_{1,\rho^*}^{\rho-1}\big)\cap\ti{U} 
=\big\{(z,v,w_k,w_+)\!\in\!\ti{U}\!:
z\!\in\!\ov\cM_{(\rho,\rho^*)}^{\rho-1};~
v\!\in\!\C^{K_{\rho^*}}\big\}.\notag
\end{gather}
It follows that
\begin{gather*}
\bigcup_{\vr^*\in f^{-1}(\rho^*)}\!\!\!\!\!\!\!\!\ov\cM_{1,\vr^*}^{\rho^+-1}
 \cap\ti{V} 
=\big\{(z,v,w_k,w_+;\ell_k)\!\in\!\ti{V}\!:
z\!\in\!\ov\cM_{(\rho,\rho^*)}^{\rho-1};\,
\ell_k\!\in\!\bP(\C^{K_{\rho^*}}\!\times\!\C)\big\} ~~\Lra\\
\bigcup_{\vr^*\in f^{-1}(\rho^*)}\!\!\!\!\!\!\!\!\ov\cM_{1,\vr^*}^{\rho^+-1} \cap\ti{U}' 
=\big\{(z,u,u_k,w_+)\!\in\!\ti{U}'\!\!: 
z\!\in\!\ov\cM_{(\rho,\rho^*)}^{\rho-1};\,
u\!\in\!\C^{K_{\rho^*}}\big\} \qquad\Lra\\
\begin{split}
& \bigcup_{\vr^*\in f^{-1}(\rho^*)}\!\!\!\!\!\!\!\!\ov\cM_{1,\vr^*}^{\rho^+}
 \cap \pi_{\rho^+,\rho^-}^{~-1}
\big(\ov\cM_{1,\rho^+}^{\rho^-}\!\cap\!\ov\cM_{1,\rho_k(j^*)}^{\rho^-}\big)
\cap\wt{W}\\
&\qquad\qquad\qquad\qquad
=\big\{(z,u,0,0;\ell)\!\in\!\wt{W}\!: z\!\in\!\ov\cM_{(\rho,\rho^*)}^{\rho-1};\,
\ell\!\in\!\bP\big(\C^{K_{\rho^*}}\!\times\!\C\big)\big\}\\
&\qquad\qquad\qquad\qquad\qquad\quad\cup 
\big\{(z,0,0,0;\ell_k)\!\in\!\wt{W}\!: z\!\in\!\ov\cM_{(\rho,\rho^*)}^{\rho-1};\,
\ell_k\!\in\!\bP\big(\C^{K_{\rho^*}}\!\times\!\C\big)\big\}.
\end{split}
\end{gather*}
Thus, by \e_ref{lift_e3} and~\e_ref{baseblow_e2}, 
\begin{equation}\label{totalspblow_e3}\begin{split}
\ti{p} \in &\big\{
f_{\rho}|_{\pi_{\rho^+,\rho^-}^{~-1}
   (\ov\cM_{1,\rho^+}^{\rho^-}\cap\ov\cM_{1,\rho_k(j^*)}^{\rho^-}) \cap \wt{W}}\big\}^{-1}
\big(\ov\cM_{1,\rho^*}^{\rho}\!\cap\!\ov\cM_{1,\rho}^{\rho}\big)\\
&\qquad\qquad\qquad
=\bigcup_{\vr^*\in f^{-1}(\rho^*)}\!\!\!\!\!\!\!\!\ov\cM_{1,\vr^*}^{\rho^+}
 \cap \pi_{\rho^+,\rho^-}^{~-1}
\big(\ov\cM_{1,\rho^+}^{\rho^-}\!\cap\!\ov\cM_{1,\rho_k(j^*)}^{\rho^-}\big)
\cap\wt{W}.
\end{split}\end{equation}\\

\noindent
In the second case above,
\begin{gather}\label{baseblow_e3}
\ov\cM_{1,\rho^*}^{\rho}\cap \ov\cM_{1,\rho}^{\rho}\cap V
=\big\{(z,0,0;\ell)\!\in\!V\!\!: 
z\!\in\!\ov\cM_{(\rho,\rho^*)}^{\rho-1};
~\ell\!\in\!\bP(\C^{K_{\rho^*}}\!\times\!0)\big\} \qquad\hbox{and}\\
\bigcup_{\vr^*\in f^{-1}(\rho^*)}\!\!\!\!\!\!\!\!\ov\cM_{1,\vr^*}^{\rho^-}\cap\ti{U} 
= f_{\rho-1}^{-1}\big(\ov\cM_{1,\rho^*}^{\rho-1}\big)\cap\ti{U}
=\ti{\cZ}_k^{\rho^-}\cup \ti{\cZ}_+^{\rho^-}, \qquad\hbox{where}\notag\\
\ti{\cZ}_{\circledast}^{\rho^-}=\big\{(z,v,w_k,w_+)\!\in\!\ti{U}\!:
z\!\in\!\ov\cM_{(\rho,\rho^*)}^{\rho-1};~
v\!\in\!\C^{K_{\rho^*}},\, w_{\circledast}\!=\!0\big\},\qquad \circledast=k,+.\notag
\end{gather}\\
We denote by $\ti{\cZ}_k^{\rho^+-1}$ and  $\ti{\cZ}_+^{\rho^+-1}$
the proper transforms of $\ti{\cZ}_k^{\rho^-}$ and $\ti{\cZ}_+^{\rho^-}$
in $\ti{V}$ and by $\ti{\cZ}_k^{\rho^+}$ and  $\ti{\cZ}_+^{\rho^+}$
the proper transforms of $\ti{\cZ}_k^{\rho^-}$ and $\ti{\cZ}_+^{\rho^-}$
in $\wt{W}$.
Then, 
\begin{gather}
\ti{\cZ}_k^{\rho^+-1}
=\big\{(z,v,0,w_+;\ell_k)\!\in\!\ti{V}\!:
z\!\in\!\ov\cM_{(\rho,\rho^*)}^{\rho-1};\,
\ell_k\!\in\!\bP(\C^{K_{\rho^*}}\!\times\!\C)\big\}\qquad\Lra \notag\\
\ti{\cZ}_k^{\rho^+-1}\cap \ti{U}' 
=\big\{(z,u,0,w_+,)\!\in\!\ti{U}'\!\!: 
z\!\in\!\ov\cM_{(\rho,\rho^*)}^{\rho-1};\,
u\!\in\!\C^{K_{\rho^*}}\big\} \qquad\Lra\notag\\
\label{totalsp_e3a}\begin{split}
\ti{\cZ}_k^{\rho^+} \cap  \pi_{\rho^+,\rho^-}^{~-1}
\big(\ov\cM_{1,\rho^+}^{\rho^-}\!\cap\!\ov\cM_{1,\rho_k(j^*)}^{\rho^-}\big)
&=\big\{(z,u,0,0;\ell)\!\in\!\wt{W}\!: z\!\in\!\ov\cM_{(\rho,\rho^*)}^{\rho-1};\,
\ell\!\in\!\bP\big(\C^{K_{\rho^*}}\!\times\!0\big)\big\}\\
&\quad\cup 
\big\{(z,0,0,0;\ell_k)\!\in\!\wt{W}\!: z\!\in\!\ov\cM_{(\rho,\rho^*)}^{\rho-1};\,
\ell_k\!\in\!\bP\big(\C^{K_{\rho^*}}\!\times\!\C\big)\big\}.
\end{split}
\end{gather}
Similarly,
\begin{gather}
\ti{\cZ}_+^{\rho^+-1}
=\big\{(z,v,w_k,0;\ell_k)\!\in\!\ti{V}\!:
z\!\in\!\ov\cM_{(\rho,\rho^*)}^{\rho-1};\,
\ell_k\!\in\!\bP(\C^{K_{\rho^*}}\!\times\!0)\big\}\quad\Lra\quad
\ti{\cZ}_+^{\rho^+-1}\cap\ti{U}'   =\eset \quad\Lra\notag\\
\label{totalsp_e3b}
\ti{\cZ}_+^{\rho^+}\cap \pi_{\rho^+,\rho^-}^{~-1}
\big(\ov\cM_{1,\rho^+}^{\rho^-}\!\cap\!\ov\cM_{1,\rho_k(j^*)}^{\rho^-}\big)
=\big\{(z,0,0,0;\ell_k)\!\in\!\wt{W}\!: z\!\in\!\ov\cM_{(\rho,\rho^*)}^{\rho-1};\,
\ell_k\!\in\!\bP\big(\C^{K_{\rho^*}}\!\times\!0\big)\big\}.
\end{gather}
Since
$$ \bigcup_{\vr^*\in f^{-1}(\rho^*)}\!\!\!\!\!\!\!\!\ov\cM_{1,\vr^*}^{\rho^+}
\cap   \pi_{\rho^+,\rho^-}^{~-1}
\big(\ov\cM_{1,\rho^+}^{\rho^-}\!\cap\!\ov\cM_{1,\rho_k(j^*)}^{\rho^-}\big)
\cap\wt{W} 
=\big(\ti{\cZ}_k^{\rho^+}\!\cap\!\ti{\cZ}_+^{\rho^+}\big)\cap 
\pi_{\rho^+,\rho^-}^{~-1}
\big(\ov\cM_{1,\rho^+}^{\rho^-}\!\cap\!\ov\cM_{1,\rho_k(j^*)}^{\rho^-}\big),$$
we conclude from \e_ref{lift_e3} and \e_ref{baseblow_e3}-\e_ref{totalsp_e3b}  
that \e_ref{totalspblow_e3} holds in this case as well.\\

\vspace{.2in}

\noindent
{\it Department of Mathematics, SUNY, Stony Brook, NY 11794-3651\\
azinger@math.sunysb.edu}\\

\end{document}